\documentclass[11pt]{amsart}

\usepackage{graphics}
\usepackage{enumitem}
\usepackage{color}
\usepackage{verbatim}

\usepackage{amsmath,amsthm, amscd, amssymb,amsfonts,mathtools}
\usepackage{wasysym}
\usepackage[all]{xy}

\usepackage[active]{srcltx}
\usepackage{array}

\usepackage{listings}

\makeindex

\newtheorem{theorem}{Theorem}[section]
\newtheorem{lemma}[theorem]{Lemma}
\newtheorem{corollary}[theorem]{Corollary}

\newtheorem{proposition}[theorem]{Proposition}
\newtheorem*{claim}{Claim}

\theoremstyle{definition}
\newtheorem{definition}[theorem]{Definition}
\newtheorem{example}[theorem]{Example}

\theoremstyle{remark}
\newtheorem{remark}[theorem]{Remark}

\newcommand\id{\operatorname{id}}
\newcommand\Id{\operatorname{id}}

\newcommand\gr{\operatorname{gr}}

\renewcommand\c{\mathfrak{c}}
\newcommand\sol{\mathfrak{s}}

\def\k{\Bbbk}

\def\ot{\otimes}

\def\s{\mathbb{S}}
\def\C{\mathcal{C}}

\def\A{\mathcal{A}}

\def\eps{\epsilon}
\def\mT{\mathcal{T}}

\def\mA{\mathcal{A}}

\def\mJ{\mathcal{J}}

\def\mH{\mathcal{H}}

\newcommand{\Up}{\Upsilon}
\newcommand{\up}{\upsilon}

\makeatletter
\newcommand*{\mathreflect}[1]{%
	\binrel@{#1}\binrel@@{\mathpalette\math@reflect{#1}}%
}
\newcommand*{\math@reflect}[2]{\reflectbox{\m@th$#1#2$}}
\newcommand*{\mathrotate}[3][]{%
	\binrel@{#3}\binrel@@{\vphantom{#3}\mathpalette\math@rotate{{#1}{#2}{#3}}}%
}
\newcommand*{\math@rotate}[2]{\math@@rotate#1#2}
\newcommand*{\math@@rotate}[4]{
	\sbox\z@{$\m@th#1#4$}%
	\smash{\makebox[\wd\z@]{\rotatebox[#2]{#3}{$\m@th#1#4$}}}%
}
\makeatother

\newcommand*{\overt}{\mathrotate[origin=c]{90}{\ominus}}

\DeclareRobustCommand\longtwoheadrightarrow
{\relbar\joinrel\twoheadrightarrow}

\DeclareRobustCommand\longrightarrow
{\relbar\joinrel\rightarrow}

\newcommand{\longhookrightarrow}{\lhook\joinrel\longrightarrow}

\newcommand{\J}{\mathbb{J}}

\newcommand{\Z}{{\mathbb Z}}

\newcommand{\N}{{\mathbb N}}

\newcommand{\R}{{\mathcal R}}

\newcommand{\Ind}{\mbox{\rm Ind\,}}

\newcommand{\Ext}{\mbox{\rm Ext\,}}

\newcommand\Rep{\operatorname{Rep}}

\newcommand{\ord}{\mathop{\rm ord}}

\renewcommand{\lg}{\langle}
\newcommand{\rg}{\rangle}

\def\pf{\begin{proof}}
	\def\epf{\end{proof}}

\newcommand\ext{\operatorname{Ext}}
\newcommand\diag{\operatorname{diag}}

\def\G{\mathbb{G}}

\def\fk{\mathcal{FK}}

\def\balpha{\boldsymbol{\alpha}}
\def\bbeta{\boldsymbol{\beta}}
\def\bdelta{\boldsymbol{\delta}}
\def\beeta{\boldsymbol{\eta}}
\def\bsigma{\boldsymbol{\sigma}}
\def\btheta{\boldsymbol{\theta}}
\def\btau{\boldsymbol{\tau}}

\def\bx{\boldsymbol{x}}
\def\by{\boldsymbol{y}}
\def\bv{\boldsymbol{v}}
\def\bw{\boldsymbol{w}}
\def\bm{\boldsymbol{m}}

\def\u{\mathfrak{u}}
\def\g{\mathfrak{g}}

\allowdisplaybreaks

\def\ben{\begin{enumerate}[leftmargin=*]}
\def\een{\end{enumerate}}

\def\bit{\begin{itemize}[leftmargin=*]}

\def\eit{\end{itemize}}

\begin{document}
	
	\title[Representations of a family of Hopf algebras]{On the representations of a family of pointed Hopf algebras}
	
	\author[Garc\'ia Iglesias, Rodriguez]{Agust\'in Garc\'ia Iglesias, Alfio Antonio Rodriguez}
	
	\address{FaMAF-CIEM (CONICET), Universidad Nacional de C\'ordoba,
		Medina A\-llen\-de s/n, Ciudad Universitaria, 5000 C\' ordoba, Rep\'ublica Argentina.}
	
	\email{agustingarcia@unc.edu.ar} 
	\email{alfio.antonio.rodriguez@mi.unc.edu.ar}
	
	\thanks{\noindent 2020 \emph{Mathematics Subject Classification.}
		16T05. \newline The work was partially supported by CONICET,
		FONCyT-ANPCyT, Secyt (UNC)}
	
	\begin{abstract}
		For each $\ell\geq 1$ and $\lambda,\mu\in\k$, we study the representations of a family of pointed Hopf algebras $\mA_{\lambda,\mu}$. These arise as Hopf cocycle deformations of the graded algebra $\fk_3\#\k\G_{3,\ell}$, where $\fk_3$ is the Fomin-Kirillov algebra and $\G_{3,\ell}$ is a given non-abelian finite group.

We compute the simple modules, their projective covers and formulate a description of tensor products. We observe that our results are fundamentally different according to the shape of the Hopf cocycle involved in the deformation.
	\end{abstract}

	\maketitle

	\section{Introduction}\label{sec:intro}
	Representation theory of non-semisimple Hopf algebras is an active area of research, as it serves a variety of interests, from Hopf algebra theory itself to application on quantum physics; meeting on its way other branches of mathematics, as modular representations of finite groups and Lie algebras. 
	
	As a distinctive, or rather defining characteristic, the category of finite dimensional representations $\Rep H$ of a Hopf algebra $H$ is a tensor category, as the coproduct determines a structure of $H$-module on the tensor product $V\ot W$, for any $V,W\in \Rep H$.
	
	As it is customary on the very reign of representation theory, the representations of a Hopf algebra $H$ can be used to have a better understanding of the algebra itself, such as its dimension or other invariants; among other characteristics, that extend up to determining if there is a Hopf cocycle relating two given Hopf algebras $H$ and $H'$ as deformations of each other. 
	
	One may in turn investigate the properties of the category $\Rep H$, starting from a description of simple modules and weight spaces, see e.g.~\cite{ARS}. As well, de Drinfeld double $D(H)$ of $H$ admits a triangular decomposition, hence some general results on such algebras can be translated to the category $\Rep D(H)$ to compute simple modules indecomposable projective modules \cite{GaV,vay-pogo}. These algebras are generally of infinite representation type; a description of the Green ring of the category $\Rep H$ of a given (tame) Hopf algebra $H$ was found in \cite{vay1}; this includes a parametrization of the indecomposable modules.  
	
	Let $\g$ be a Lie algebra and let $U_q(\g)$, resp.~$u_q(\g)$ denote the quantum group, resp.~small quantum group, given by the quantum enveloping algebra of $\g$; here $q$ is a root of 1. Then $\Rep U_q(\g)$ is governed by the symplectic leaves of the Poisson-Lie groups $G^d$ (connected) and $G$; for $\operatorname{Lie} G^d \simeq  \g^\ast$, \cite{dkp, dl, ls}.
	As well, the category $\Rep \u_q (\g)$, with $\ord q = p$ prime is analogous to $\Rep u(\g)$, in char. $p > 0$, \cite{l2}. For a generic parameter $q$, the category $\Rep U_q(\g)$ is closely related to $\Rep U(\g)$ \cite{l1, rosso}.
	
	\subsection{Main results}
	
	Fix $\ell\in\N$. In this paper we look into the representation theory for a family of finite-dimensional pointed Hopf algebras $\mA_{\lambda,\mu}$, see \eqref{eqn:definitionA}, over a group $\G_{3,\ell}$, $\lambda,\mu\in\k$. We list their simple modules in Theorem \ref{thm:simple}, together with the corresponding projective covers in Theorem  \ref{thm:projective}. We show that these algebras are of infinite representation type. 
	
	When $\ell=1$, we get that $\G_{3,1}\simeq \s_3$ and these results replicate some contributions from \cite{GI} for the finite-dimensional pointed Hopf algebras $\mH_\lambda$, $\lambda\in \k$, over $\s_3$. These algebras can be described as $\mH_\lambda\simeq \mA_{\lambda,0}$ in this context; hence our results apply. 
	
	As a matter of fact, we show that the tensor category $\R=\Rep\mA_{\lambda,\mu}$ is graded and the grading $\R=\bigoplus_{j=0}^{\ell-1} \R_j$ is such that $\R_0$ is tensor equivalent to $\Rep \mH_\lambda$.
	
	This is the content of Proposition \ref{pro:tensor} and Corollary \ref{cor:tensor}, where the tensor structure of certain family of $\mA_{\lambda,\mu}$-modules -which includes, but not describes, the indecomposable ones- is set.
	
	We summarize the classification results in  Theorem \ref{thm:simple} facts, avoiding both {\it ad-hoc} notation from the aforementioned theorems and some particular sub-cases, in the following table:
\begin{center}
	\begin{tabular}{| c | c | c|}
	\hline
&$\mu={\lambda}/{3}$ & $\mu\neq {\lambda}/{3}$ \\ 
	\hline
		dimension & \multicolumn{2}{c|}{number} \\ 
\hline
	1 & $2$ &  $2$ \\ 
	\hline
	2 & $5\ell-1$  & $4$ \\ 
	\hline
	6 & $-$  & $2\ell-2$ \\ 
\hline
\end{tabular}
\end{center}
As for the exceptions, in the column $\mu\neq\lambda/3$, we mention that:
\begin{itemize}[leftmargin=*]
\item When $\lambda=0$, there is a unique two-dimensional module when $\lambda=0$, as the four modules in the table collapse into a single one. 
\item As well, for each $\ell$ there is a collection $\{\c_{\pm j}:j=1,\dots,\ell-1\}\in \k\setminus \{1\}$, see \eqref{eqn:c}, of exceptional scalars  for which the case $\mu=\c_{\pm j}\lambda/3$ makes two (and only two) 6-dimensional modules collapse into one, and hence the number of such modules is $2\ell-3$.
\end{itemize} 
We refer to Section \ref{sec:classification}, particularly to our main result Theorem \ref{thm:simple} for further detail.
%
	
\subsection{Our moonshine}
	
		The line $\mu=\frac{\lambda}{3}$  is precisely the space in which the Hopf cocycle involved in the deformation of the graded algebra $\mA_{0,0}$ into $\mA_{\lambda,\mu}$ can be obtained as an exponential of a Hochschild 2-cocycle on $\mA_{0,0}$. Outside this line, the cocycles involved are pure. 
	
		As mentioned above, this line also splits the structure of the category $\R$, from the very number and shape of simple modules. We do not wish to dwell on these concepts here, and we refer to \cite{GS} for context; we simply wish to point out the occurrence of this relation- in these two seemingly unrelated (or a priori hardly related) contexts: our own little moonshine.
	
\subsection{General structure}	
The algebras $\mA_{\lambda,\mu}$- more generally any (known) finite-dimensional pointed Hopf algebra, naturally posses the following characteristics:
	\begin{itemize}
		\item There is a group $G$ so that $\k G$ is a subalgebra of $A$.
		\item There is a set $X$ of elements in $A$ and a set of monic $\N_0$-homogeneous words $L$ on $X$ so that $\{ag:g\in G, a\in L\}$ is a linear basis of $A$. 
	\end{itemize}
	Consider the decomposition of the left regular representation ${}_G\k G$ of $G$ by the irreducible representations $\widehat{G}$, namely $\k G\simeq  \bigoplus_{S\in\widehat{G}}S^{d_S}$, where $d_S=\dim_\k S$.
	In particular, the left regular representation ${}_AA$ of $A$ can be decomposed as
	\[
	{}_AA \simeq {}_AA_G\ot_G {}_G\k G  \simeq {}_AA_G\ot_G \bigoplus_{S\in\widehat{G}}S^{d_S}\simeq \bigoplus_{S\in\widehat{G}}(A\ot_G S)^{d_S},
	\]
	where ${}_AA_G$ denotes the $(A,\k G)$-bimodule $A$ defined by left and right multiplication, respectively. Hence the induced modules $A(S)\coloneqq A\ot_G S$ are straightforward examples of projective $A$-modules. 
	
	We remark that these are not necessarily projective. If $B_S$ is a basis of the $G$-module $S$, then $\{av:a\in L,v\in B_S\}$ is a basis for $A(S)$. The indecomposable projective modules for $A$ are in one-to-one correspondence with the irreducible modules for $A$ via:
	\begin{align*}
		L&\mapsto P(L)\text{ projective cover of }L, & P&\mapsto \text{irreducible quotient } P\twoheadrightarrow L_P.
	\end{align*} 
	Moreover, the following decomposition holds:
	\[
	{}_AA\simeq \bigoplus\limits_{L\text{ irreducible}}P(L)^{\dim L}.
	\]
	This is very useful, in a numerological sense, to determine the number of simple modules to prove that a given module is a projective cover:
	\begin{align}\label{eqn:dim-simple-proy}
		\dim_\k A=\sum_{L}\dim_\k L\dim_\k P(L).
	\end{align}
	We are interested in finding the simple $A$-modules and their projective covers. We propose the following strategies, based both on our experience in the subject and some standard considerations:
	\begin{enumerate}[leftmargin=*]
		\item If $S\in\widehat{G}$, then check if there is an action of $A$ on $S$ compatible with the $G$-structure. In this case, $A(S)$ projects onto $S$. If $A(S)$ is indecomposable, then it is the projective cover of $S$ as $A$-module \cite[Proposition 5.4 (ii)]{GI}. Otherwise, it can be computed as a factor of $A(S)$. 
		\item Study restrictions on sums $S_1\oplus\dots\oplus S_k$, of (non-necessarily) modules $S_i\in\widehat{G}$, to encompass an action of $A$ extending, once again, that of $G$.
		\item Find the $G$-module structure of $A(S)$, for each $S\in\widehat{G}$, and analyze:
		\begin{enumerate}
			\item wether this module is indecomposable, or even simple,
			\item irreducible quotients $A(S)\twoheadrightarrow L$.
		\end{enumerate} 
	\end{enumerate}
	It goes without saying that this includes a step zero: 
	\begin{itemize}[leftmargin=1cm]
		\item[(0)] Find the irreducible representations of $G$.
	\end{itemize}
	Moreover, any irreducible $A$-module $L$ arises as one of the quotients described in 3b above. Indeed, if $L$ is simple and $L\simeq S_1\oplus \dots \oplus S_k$ as $G$-module, then the projection $A\ot_GL=\bigoplus_{i=1}^k A(S_i)\twoheadrightarrow L$ gives rise to morphisms $A(S_i)\hookrightarrow  A\ot_GL \twoheadrightarrow L$ which cannot be all zero and thus determine (at least) an index $j$ with $A(S_j)\twoheadrightarrow L$. 
	
\subsection{Organization}
The paper is organized as follows. In Section \ref{sec:preliminaries} we lay out some fundamental concepts and basic results on representation theory that will be used in the sequel. As well, we describe the algebras $\mA_{\lambda,\mu}$ and review the main characteristics of the Hopf cocycles involved. As well, we compute the simple modules for the group(s) $\G_{3,\ell}$, together with their duals and tensor products.

We compute the simple $\mA_{\lambda,\mu}$-modules in Section \ref{sec:simples}; we write our classification in \S\ref{sec:classification} and describe the tensor structure of this category in \S\ref{sec:tensor-structure}. 

Finally, we compute the projective covers in Section \ref{sec:projective} and study the representation type of $\mA_{\lambda,\mu}$ in Section \ref{sec:extensionsdiagrams}, where we compute extensions between simple modules in \S\ref{sec:extensions} and obtain the corresponding Gabriel quivers in \S \ref{sec:diagrams}.

\section{Preliminaries}\label{sec:preliminaries}

We work over an algebraically closed field $\k$ of characteristic zero.

\subsection{Representation of algebras}
Let $A$ be a finite-dimensional $\k$-algebra and $M$ be a left $A$-module. A \emph{projective cover} of
$M$ is a pair $(P(M),f)$ with $P=P(M)$ a projective $A$-module and
$f:P\to M$ an\emph{ essential map}: $f$ is surjective and
for every proper submodule
$N\subset M$ it follows that $f(N)\neq M$. Projective covers
are unique up to isomorphism and always exist if $\dim A<\infty$. See \cite[Sect.
6]{CR} for details. A key feature of projective covers is formula \eqref{eqn:dim-simple-proy}, which is a helpful tool to determine the number of simple modules. 

In turn, to study the representation type of $A$, we look at the spaces $\Ext^1(S,S')$ of extensions between each pair of simple modules of $A$. We briefly review this tool. Recall that  $\Ext^1(S,S')$ is the set of (equivalence classes of) modules $M$ fitting in an exact sequence $0\to S'\to M\to S\to 0$. The Baer sum endows this set with a $\k$-linear structure (so that $0=S\oplus S'$); in particular we can compute $d_{S,S'}=\dim \Ext(S,S')\in\N\cup\{0\}$. The {\it Gabriel Quiver} $G=G(A)$ associated to $A$ is the directed graph set with of vertices $V(G)$ given by one vertex per isomorphism class $[S]$ of simple $A$-modules and $d_{S,S'}$  arrows $[S]\to [S']$. If $Q$ is any quiver, with vertices $V=V(Q)=\{1,\dots,n\}$,  the {\it separated graph} $S=S(Q)$ of $Q$ is the unoriented graph obtained by duplicating the set of vertices of $G$, so $V(S)=V'\cup V''=\{1',1'',\dots,n',n''\}$, for $V'=V''=V$ and with an edge $i'-j''$ for each arrow $i\to j$ in $Q$. The separated graph of an algebra $A$ is the separated graph of $G(A)$.

By \cite[Theorem 2.6]{AuRS}, it follows that $A$ is
of finite (tame) representation type if and only if its separated
diagram is a disjoint union of finite (affine) Dynkin
diagrams.

We recall from \cite[\S 4.14]{egno} that a {\it grading} of a tensor category $\mathcal{C}$ by a group $G$ is a decomposition $\C=\bigoplus_{g\in G}$ into a direct sum of abelian subcategories  $\C_g\subset \C$, $g\in G$, and such that $\ot$  maps $(\C_g, \C_h)$ into  $C_{gh}$. If $e\in G$ denotes the unit, then the subcategory $\C_e$ is a tensor subcategory, called the
trivial component of the grading. In this setting $\C$ is called an extension of $\C_e$. The grading is said to be faithful if $\C_g\neq 0$ for all $g \in G$.

\subsection{Pointed Hopf algebras over $\G_{3,\ell}$ -and Hopf cocycles}
For each $\ell\geq 1$, we consider the group
\begin{align}\label{eqn:group-def}
	\mathbb{G}_{3,\ell}&=\langle s,t|s^3=t^{2\ell}=1, ts=s^2t\rangle=C_3\rtimes C_{2\ell},
\end{align}
where $C_3=\langle {s} \rangle, C_{2\ell}=\langle t \rangle$ are cyclic groups of orders $3$ and $\ell$, respectively. 
Observe that $|\mathbb{G}_{3,\ell}|=6\ell$.
If $\ell=1$, then $\mathbb{G}_{3,\ell}\simeq \s_3$.

\noindent We fix the {\it bosonization}
	$\mT=\k\langle a_0, a_1, a_2\rangle \# \k \mathbb{G}_{3,\ell}$,
where $\#$ is determined by
\begin{align}\label{eqn:commutation}
	t\cdot a_i&=-a_{-i}, && 	s\cdot a_i=a_{i+1}, \qquad i\in X\coloneqq \{0,1,2\}.
\end{align}
Here, we write $a_i$ instead of $a_i\#1$ and $t,s$ for $1\#t$, $1\#s$, for a clean exposition. As well, we identify  $X\leftrightarrow\Z/3\Z=\{\bar{0},\bar{1},\bar{2}\}$.

Next, we introduce the ideal $\mJ_{\lambda,\mu}\subset \mT$ generated by
\begin{align}
	\notag	a_0^2=a_1^2=a_2^2&=\mu(1-t^2),\\
	\label{eqn:rels_def}	a_0a_1+a_1a_2+a_2a_0&=\lambda(1-st^2),\\
	\notag 	a_1a_0+a_2a_1+a_0a_2&=\lambda(1-s^2t^2).
\end{align}
We shall consider the family of algebras $\A_{\lambda,\mu}$ defined as the quotient
\begin{align}\label{eqn:definitionA}
	\A_{\lambda,\mu}=\k\langle a_0, a_1, a_2\rangle \# \k \mathbb{G}_{3,\ell} /\mJ_{\lambda,\mu}. 
\end{align}
These are pointed Hopf algebras with coradical $\k \mathbb{G}_{3,\ell}$ and comultiplication 
\[
\Delta(a_i)=a_i\ot 1 + s^{-i}t\ot a_i, \ i\in X.
\]
In particular, the antipode $S\colon \A_{\lambda,\mu}\to \A_{\lambda,\mu}$ is determined by the formulae 
\begin{align*}
S(g)&=g^{-1},  \ g\in\G_{3,\ell}, & S(a_i)&=a_{i}t^{-1}s^i, \ i=0,1,2.
\end{align*}

In order to deal with representations of $\A_{\lambda,\mu}$, it is important to keep in mind that, as an algebra, it is generated by $\{a_0,s,t\}$ by \eqref{eqn:commutation}. Hence the task can be thought of finding three matrices $[a_0]$, $[s]$, $[t]$ satisfying some identities. Namely, those involving $s$ and $t$ from \eqref{eqn:group-def}, together with \eqref{eqn:rels_def}.

\begin{remark}\label{rem:ell=1}
We shall generically assume that $\ell>1$, as otherwise $\mathbb{G}_{3,\ell}\simeq \s_3$ if $\ell=1$ and the representation theory of the pointed Hopf algebras over this group has been studied in \cite{GI}. Nevertheless, the results on this article encompass those in loc.cit.
\end{remark}

We collect a few of important remarks about these algebras:
\begin{itemize}[leftmargin=*]
\item For any $k=0,\dots,\ell-1$ a Hopf algebra $\A_{\lambda,\mu}^{(k)}$ can be defined, in such a way that 
$\A_{\lambda,\mu}=\A_{\lambda,\mu}^{(0)}$. Indeed, $\A_{\lambda,\mu}^{(k)}=\A_{\lambda,\mu}$ as algebras but the comultiplication is now given by 
\[
\Delta(a_i)=a_i\ot 1 + s^{-i}t^{2k+1}\ot a_i, \ i\in X.
\]
This does not affect the representation theory of the algebras -as far as abelian categories are concerned. It does impact the tensor structure, which we expect to address in future work. 

\item The family $\A_{\lambda,\mu}^{(k)}$ is a complete list of finite-dimensional pointed Hopf algebras over $\mathbb{G}_{3,\ell}$ and such that their infinitesimal braiding is a principal realization of the affine rack $\operatorname{Aff}(3,-1)$, see \cite{GV}.

\item The graded object $\gr \A_{\lambda,\mu}^{(k)}$ associated to the coradical filtration of these algebras is the bosonization $\fk_3\#\k\G_{3,\ell}$ of the Fomin-Kirillov algebra $\fk_3$ with the group algebra $\k\G_{3,\ell}$. 

Recall that $\fk_3=\k\lg x_0,x_1,x_2\rg$ with relations:
\begin{align*}
x_0^2&=x_1^2=x_2^2=0, & x_0x_1+x_1x_2+x_2x_0&=x_1x_0+x_2x_1+x_0x_2=0.
\end{align*} 

\item $\A_{\lambda,\mu}^{(k)}$ is a Hopf cocycle deformation of $\fk_3\#\k\G_{3,\ell}$. The corresponding cocycle can be obtained as an exponential $e^\eta$ of a Hochschild 2-cocycle $\eta$ if and only if $\mu=\frac{\lambda}{3}$. Otherwise the cocycle is so-called pure and can be computed by ad-hoc combinatorics. We refer to \cite{GS} for details.
\end{itemize}

\subsection{Irreducible representations of $\G_{3,\ell}$}

We fix $\ell> 1$ and set $G=\G_{3,\ell}$. Let us fix as well $\xi$ and $\zeta\in \mathbb{C}$ which are primitive roots of one of order 3 and $\ell$, respectively.
For each $j= 0, \dots, \ell-1$ we fix $b_j\in\k$ with $b_j^2=\zeta^j$.

Next proposition we describe the set of irreducible representations of $G$, see \S \ref{sec:proofG} for the proof.
\begin{proposition}\label{pro:simplesG}
The irreducible representations of $\mathbb{G}_{3,\ell}$ are given, up to isomorphism by the following modules.
	\begin{itemize}[leftmargin=*]
		\item The one-dimensional modules  $S_{j}^\pm=\lg z_{\pm}\rg$, $j= 0, \dots, \ell-1$, with action
		\begin{align}\label{eqn:accion-Sj}
			s\cdot z_{\pm} &= z_{\pm}, & t\cdot z_{\pm} &= \pm b_j z_{\pm}.
		\end{align}
		\item The two-dimensional modules $M_{j}=\lg v,w\rg$, $j= 0, \dots, \ell-1$, with action 
		\begin{align}\label{eqn:accion-Mj}
			[s]:=\begin{pmatrix*} \xi&0\\0&\xi^{2} \end{pmatrix*}, \qquad [t]:=\begin{pmatrix*} 0&\zeta^{j}\\ 1&0 \end{pmatrix*}.
		\end{align}
	\end{itemize}
\end{proposition}
We define the simple modules of {\it type $j$} as the subset
\[
\widehat{G}[j]=\{S_j^+,S_j^-,M_j\}\subset \widehat{G},
\]
so $\widehat{G}=\bigcup_{j=0}^{\ell-1}\widehat{G}[j]$. 
As well, for a given $G$-module $N$ we set, for each $j$, the submodule generated by the irreducible submodules of type $j$:
\begin{align}\label{eqn:j-part}
	N[j]\coloneqq\lg S\subseteq N : S\in \widehat{G}[j]\rg.
\end{align}
In particular $N=\bigoplus\limits_{j=0}^{\ell-1}N[j]$. The following observation will be very useful.
\begin{lemma}\label{lem:uniquej}
Let $M$ be a $G$-module. Then 
\begin{align*}
x\in M[j] \text{ if and only if }t^2\cdot x=\zeta^j\,x.
\end{align*}
\end{lemma}
\pf
One implication is clear. For the other, let $x\in M$ be such that $t^2\cdot x=\zeta^j\,x$ and write $x=\sum_{i=0}^{\ell-1}x_i$, with $x_i\in M[i]$.
Then 
\[
\sum_{i=0}^{\ell-1}\zeta^j\,x_i=\zeta^j\,x=t^2\cdot x=\sum_{i=0}^{\ell-1}\zeta^i\,x_i,
\]
which shows that $x_i=0$ if $i\neq j$ and thus $x=x_j\in M[j]$.
\epf

\subsubsection{Computations of $G$-modules}

In this part we construct a family of $G$-modules, in order to pin down the irreducible representations in $\widehat{G}$.
We denote by $C_\theta$ the cyclic group of order $\theta\in\N$.

\begin{lemma}
The group $G$ fits in the exact sequence
\[1\longrightarrow H \longrightarrow G \longrightarrow C_2 \longrightarrow 1.\]
where $H\simeq C_3\times C_\ell$ if $3\mid \ell$ and $H\simeq C_{3\ell}$ if $3\nmid \ell$.
\end{lemma}
\pf
Let $u=t^2\in G$, so we can write 
\begin{align*}
	G=\langle s, u, t| su=us, s^3=u^\ell=1, t^2=u, tu=ut, ts=s^2t \rangle.
\end{align*}
and the lemma follows for $H=\langle s, u \rangle \lhd G$.
\epf

As $H\coloneqq \lg s,u\rg\leq G$ is an abelian group, we have that
\begin{align*}
	\widehat{H}= \{\sigma_i\times\rho_j\,|\, i=0, 1, 2, \ j=1, \dots, \ell-1\},
\end{align*}
for $\sigma_i\times\rho_j(s)=\xi^{i}$, $\sigma_i\times\rho_j(u)=\zeta^{j}$. Let us set $\lg x\rg$ to denote the 1-dimensional $H$-module associated to a given representation  $\sigma_i\times\rho_j\in  \widehat{H}$.
Now
\begin{align*}
	\Ind_H^G(\sigma_i\times\rho_j)&=G\ot_H \lg x\rg=\lg 1\ot_H  x,  t\ot_H x\rg,	
\end{align*}
and thus, if $v\coloneqq 1\ot x$ and $w\coloneqq t\ot x$, we have that
\begin{align*}
	s\cdot(1\ot x) &= \xi^i (1\ot x)= \xi^i v, && t\cdot (1\ot x) = t\ot x= w,\\
	s\cdot (t\ot x) &= \xi^{2i} (t\ot x)= \xi^{2i} w, && t\cdot (t\ot x) = \zeta^{j} (1\ot x)= \zeta^{j} v.
\end{align*}

\begin{lemma}
Let $M_{i,j}=\lg v,w\rg$ be the $G$-module with action determined by
	$[s]:=\begin{psmallmatrix*} \xi^{i}&0\\0&\xi^{2i} \end{psmallmatrix*}$, $[t]:=\begin{psmallmatrix*} 0&\zeta^{j}\\ 1&0 \end{psmallmatrix*}$.
Then:
\begin{enumerate}[leftmargin=*]
\item $M_{0,j}=S_j^+\oplus S_j^-$, where $S_j^\pm$ are the simple $G$-modules generated by
\[
z_+=b_jv+w \qquad \text{and } \qquad z_-=-b_jv+w
\]
respectively. In particular, $s\cdot z_{\pm}=z_{\pm}$ and $t\cdot z_{\pm}=\pm b_j\,z$.
\item Set $M_{j}\coloneqq M_{1,j}$. Then $M_j$ is irreducible and $M_{2,j}\simeq M_{j}$.
\end{enumerate} 
\end{lemma}
The modules $S_j^\pm$ and $M_j$ coincide with those in Proposition \refeq{pro:simplesG}.
\pf
{\it (1)} We look for a one-dimensional submodule of $M_{i,j}$, set 
$z=\alpha v + \beta w$ and assume that $s\cdot z=a\,z$, $t\cdot z=b\,z$, for some  $a,b \in \k$. Observe that $\alpha\beta\neq 0$, as the action of $t$ interchanges $v$ and $w$. Now we get
\begin{align*}
a\,z&=a\alpha v+a\beta w= s\cdot z = s\cdot(\alpha v + \beta w) = \alpha\xi^i v + \beta \xi^{2i} w,\\
b\,z&= b\alpha v+b\beta w=t\cdot z = t\cdot(\alpha v + \beta w) =  \beta \zeta^{j} v +\alpha w.
\end{align*}
Then $	\alpha\xi^i = a\alpha \Rightarrow a=\xi^i$ and $\beta\xi^{2i} = a\beta \Rightarrow a= \xi^{2i}$,
and we get that we need $i=0$ for this, and hence $a=1$. Now, we also have	
$
\alpha = b\beta$ and 
$\beta\zeta^{j} = b\alpha$, 
which implies $b^2=\zeta^j$, $\alpha=b\beta$. Fix $b=b_j$ so that $b_j^2=\zeta^j$. Then
$z_+=b_jv+w$ and $z_-=-b_jv+w$ define submodules $S_j^+=\lg z_+\rg$ and $S_j^-=\lg z_- \rg$ as in the statement of the lemma.

{\it (2)} We have just seen above that these modules are irreducible, as they do not admit one dimensional submodules.
Let us denote by $B_i=\{v_i, w_i\}$ the basis of $M_{i,j}$, $i=1,2$. Then 
$v_1\mapsto w_2$, $w_1\mapsto \zeta^{j} v_2$ defines and $G$-module isomorphism $M_{1,j}\to M_{2,j}$.
\epf

\subsubsection{Proof of Proposition \ref{pro:simplesG}}\label{sec:proofG}
We start by observing that
\[
|\G_{3,\ell}|=6\ell=2\ell \underbrace{(\dim S_{j}^\pm)^2}_{=1} + \ell \underbrace{(\dim M_{j})^2}_{=4}.
\]
Hence the result follows as $S_j^+\not\simeq S_j^-$ and $S_j^+\not\simeq S_k^\pm$, $k\neq j$.
\qed

\subsubsection{On the tensor structure of $\Rep \G_{3,\ell}$}
Next result will be useful on \S\ref{sec:tensor-structure}.

\begin{proposition}\label{pro:tensor-group} The following identities hold.
	\begin{enumerate}
		\item $(S_j^{\pm})^\ast\simeq S_{\ell-j}^\pm$ and $M_j^\ast\simeq M_{\ell-j}$.
		\item $S_j^\eps\ot S_k^{\eps'}\simeq S_{j+k}^{\eps\eps'}$ and $S_j^-\ot M_k\simeq M_k\ot S_j^-\simeq M_{j+k}$.
		\item $M_j\ot M_k\simeq S_{j+k}^+\oplus S_{j+k}^-\oplus M_{j+k}$.
	\end{enumerate}
\end{proposition}
\pf
Set $M_j=\lg v,w\rg$. If $\{f,g\}$ is the dual basis of $\{v,w\}$, then it follows that $\lg f,g\rg \simeq M_{\ell-j}$. This shows {\it (1)} The first part of {\it (2)} is clear, for the second, if $S_j^-=\lg y\rg$ and $M_k=\lg u,z\rg$, then $S_j^- \ot M_k =\lg y\ot u, -b_jy\ot z\rg\simeq M_{j+k}$.
As for {\it (3)}, we see that $\lg b_j\,v\ot z\pm b_k w\ot u\rg\simeq S_{j+k}^\pm$ and $\lg w\ot z,v\ot u\rg\simeq M_{j+k}$.
\epf

\subsection{The graded case}
We shall assume along the work that $\lambda$ and $\mu$ are not simultaneously zero, as this leads to the graded case $\A_{0,0}=\fk_3\#\k\G_{3,\ell}$ which lands in the general picture described in \cite[Proposition 4.1]{GI}, that in this context reads as follows.
\begin{proposition}\label{pro:graded-case}
	The simple modules $L$ for $\mA=\mA_{0,0}$ are in bijective correspondence $L=L_S\leftrightarrow S$
	with the simple modules $S$ over $\G_{3,\ell}$.

Given $S\in\widehat \G_{3,\ell}$, $L_S$ is the
	$\mA$-module such that
	$$
	L_S\cong S \text{ as } \G_{3,\ell}\text{-modules,} \quad\text{and  }\quad  a_0\cdot L_S=0.
	$$
	This correspondence preserves tensor products and duals.
\end{proposition}

\subsection{Notation}\label{sec:notation}
Whenever $\ell>1$ is fixed, we may write $G$ to refer to $\G_{3,\ell}$. 
As well, we let  $\xi$ and $\zeta$ denote primitive roots of $1$ of orders $3$ and $\ell$, respectively. 
We set $\J=\{1,\dots,\ell-1\}$, and $\J_0=\J\cup\{0\}$.
Finally, we fix a choice of scalars $b_j\in\k$,  one for each $j\in\J_0$, so that $b_j^2=\zeta^j$; with $b_0=1$.  

If $M$ is an $\mA_{\lambda,\mu}$-module, then we shall write $M_{|}$ to denote its restriction as a representation of $G$. If $S\in\widehat{G}$ and $M_{|}\simeq S$, then we refer to such $\mA_{\lambda,\mu}$-module $M$ as an {\it extension} of $S$.

\section{Simple $\mA_{\lambda,\mu}$-modules}\label{sec:simples}

Fix $\ell>1$. In this section we classify, up to isomorphism, the simple modules for the algebras $\mA_{\lambda,\mu}$. We start by determining in \S\ref{sec:simplesoverG} those simple $\A_{\lambda, \mu}$-modules $L$ such that their restriction $L_{|}$ are simple modules for the group $G$.  Next, we study general modules in \S\ref{sec:mixed}.
We  write our classification results in \S\ref{sec:classification}.
We consider the notation from \S\ref{sec:notation}.

\subsection{Simple $\mA_{\lambda,\mu}$-modules defined over $\widehat{G}$}\label{sec:simplesoverG}

We study extensions of 1-dimensional $G$-modules in \S\ref{sec:ext-simples1} and of the 2-dimensional ones in \S\ref{sec:ext-simples2}.

\subsubsection{Extensions of the one-dimensional modules $S_j^\pm$.}\label{sec:ext-simples1}

We study modules $L$ so that $L_{|}\simeq S_j^\pm$. We show that necessarily $j=0$ and both $S_0^+$ and $S_0^-$ admit a (unique) extension.

\begin{lemma}\label{lem:nontrivial action dim1}  
The trivial action $a_0=0$ defines a structure of $\mA_{\lambda,\mu}$-module $L_0^\pm$ extending the $G$-modules $S_0^\pm$.

Conversely, if  $L$ is a simple  $\A_{\lambda, \mu}$-module of dimension 1 then $L\simeq L_0^\pm$.
\end{lemma} 
\pf
The first part is clear. Let $L_{|}=\lg w\rg\simeq S_j^\pm$, so that $s\cdot w=w$ and $t\cdot w=b_jw$.
Let $\alpha\in\k$ be such that $a_0\cdot w=\alpha w$. Since $ta_0=-a_0t$, this rapidly gives $\alpha=0$. Thus 
$0=	a_0^2\cdot w=\mu(1-t^2)\cdot w=\mu(1-\zeta^j)w$
and therefore $\mu=0$ or $j=0$. As $a_1=sa_0s^{2}$, $a_2=s^{2}a_0s$, we have
\[
0=(a_0a_1+a_1a_2+a_2a_0)\cdot w= \lambda(1-st^2)\cdot w=\lambda(1-\zeta^j)w
\]
so that $\lambda=0$ or $j=0$.
As $(\lambda,\mu)\neq(0,0)$, we get  $j=0$ and thus $L\simeq L_0^\pm$. 
\epf

\subsubsection{Extensions of the two-dimensional modules $M_j$.}\label{sec:ext-simples2}

Next we look at  modules $L$ with $L_{|}\simeq M_j$. We consider such $G$-modules together with a fixed basis $\{v,w\}$ so that \eqref{eqn:accion-Mj} holds.

For each couple $a,c\in \{\sqrt{\lambda/3},-\sqrt{\lambda/3}\}$ we consider the matrices
\begin{align}\label{eqn:matrixmj}
m_j(a,c)\coloneqq\begin{pmatrix*} a&-\zeta^jc\\c&-a \end{pmatrix*}, \qquad j\in\J_0.
\end{align}
We include $m_j(0,0)=\begin{psmallmatrix}
	0&0\\0&0 \end{psmallmatrix}$ for $\lambda=0$. As well, $m_0(a,c)=\begin{psmallmatrix}
	a&-c\\c&-a \end{psmallmatrix}$.

\begin{lemma}\label{lem:LMj}
Fix $a,c\in\{\pm\sqrt{\lambda/3}\}$.
\begin{enumerate}[leftmargin=*]
\item If $j=0$, then there is an $\mA_{\lambda,\mu}$-module $L=L_0(a,c)$ such that $L_{|}\simeq M_0$ with action, on the basis $\{v,w\}$ given by the matrix $m_0(a,c)$.

\item If $j\in\J$ and $\mu=\frac{\lambda}{3}$, then there is an $\mA_{\lambda,\mu}$-module $L=L_j(a,c)$ such that $L_{|}\simeq M_j$ with action, on the basis $\{v,w\}$ given by the matrix $m_j(a,c)$.
\end{enumerate}
Conversely, if $L$ is an $\mA_{\lambda,\mu}$-module such that $L_{|}\simeq M_j$, then $L\simeq L_j(a,c)$ for some $a,c\in\{\pm\sqrt{\lambda/3}\}$.

In any case $L_j(a,c)\simeq L_j(a',c')$ if and only if  $(a,c)=(a',c')$.
\end{lemma}
\pf
It is straightforward to check that, for the conditions in the statement, \textit{(1)} and \textit{(2)} hold. The same applies for isomorphism classes.

For the converse, let us set $[a_0]=\begin{pmatrix*} a&b\\c&d \end{pmatrix*}$ to represent the action of $a_0$. Now, as $ta_0\stackrel{!}{=}-a_0t$ we get that $a=-d$ and $\zeta^jc=-b$. Indeed:
\[
\begin{pmatrix*} 0&\zeta^j\\1&0\end{pmatrix*}\begin{pmatrix*} a&b\\c&d\end{pmatrix*}=\begin{pmatrix*} \zeta^jc&\zeta^jd\\a&b\end{pmatrix*}\stackrel{!}{=}-\begin{pmatrix*} b&\zeta^ja\\d&\zeta^jc\end{pmatrix*}=-\begin{pmatrix*} 0&\zeta^j\\1&0\end{pmatrix*}\begin{pmatrix*} a&b\\c&d\end{pmatrix*}\]
Now, $[a_0]=\begin{psmallmatrix*} a&-\zeta^jc\\c&-a\end{psmallmatrix*}\implies [a_0]^2=\begin{psmallmatrix*} a^2-\zeta^jc^2&0\\0&a^2-\zeta^jc^2\end{psmallmatrix*}$.
As this must coincide with $\mu(\id-[t]^2)=\mu(1-\zeta^j)\id$, we have that
\begin{align}\label{eqn:dim2-cond1}
	a^2-\zeta^jc^2=\mu(1-\zeta^j).
\end{align}
Now, 
$
[a_1]=[s][a_0][s]^{-1}=\begin{psmallmatrix*} a&-\zeta^j\xi^2 c\\ c\xi &-a \end{psmallmatrix*}$,
$[a_2]=[s]^{-1}[a_0][s]=\begin{psmallmatrix*} a&-\zeta^j\xi c\\ c\xi^2 &-a \end{psmallmatrix*}$, so
\[
[a_0][a_1]+[a_1][a_2]+[a_2][a_0]=\begin{pmatrix*}3(a^2-\xi c^2\zeta^j)&0\\0&3(a^2-\xi^2 c^2 \zeta^j)\end{pmatrix*}\]
which must coincide with the matrix $\lambda(1-[s][t]^2)=\begin{psmallmatrix*} \lambda(1-\xi\zeta^j)&0\\0&\lambda(1-\xi^2\zeta^j)\end{psmallmatrix*}$.
Thus $a^2=c^2=\frac{\lambda}{3}$, which combined with \eqref{eqn:dim2-cond1} yields:
$\mu(1-\zeta^j)=\frac{\lambda}{3}(1-\zeta^j)$,
and therefore for the existence of such a module we necessarily have $\mu = \frac{\lambda}{3}$ or $\zeta^j=1$, namely $j=0$. Hence $M\simeq L_j(a,c)$.
\epf

\subsection{``Mixed'' $\mA_{\lambda,\mu}$-modules}\label{sec:mixed}

In this section we study simple modules $L$ for $A_{\lambda,\mu}$ when $L_{|}\notin \widehat{G}$. 
We proceed as follows: In \S\ref{subsec:uniquej} we show that if $L$ is indecomposable then $L_{|}=L_{|}[j]$ for some $j\in\J_0$, see \eqref{eqn:j-part}. Next, in \S\ref{subsec:powersof1} we study powers of the one dimensional modules, that is when $L_{|}\simeq S_j^{p}\oplus S_j^{-q}$, $p,q\geq 0$. In \S\ref{subsec:powersof2} we investigate powers of the two dimensional module, that is the case $L_{|}\simeq M_j^{r}$, $r\geq 0$. Finally in \S\ref{subsec:verymixed} we look at the general ``very mixed'' case $L_{|}\simeq S_j^{p}\oplus S_j^{-q}\oplus M_j^{r}$, $p,q,r\geq 0$.

 In \S\ref{subsec:verymixed6} we describe a 6-dimensional mixed irreducible module $\mA$-module  such that $L_{|}\simeq S_j^+\oplus S_j^{-}\oplus M_j^{2}$.	

\subsubsection{A single $j$}\label{subsec:uniquej}

We recall the definition \eqref{eqn:j-part} of $N[j]$ for a $G$-module $N$, which determines $N=\bigoplus_{j=0}^{\ell-1}N[j]$. In particular, any $S\in\widehat{G}$ is such that $S=S[j]$ for some $j\in\J_0$.
For a $\mA_{\lambda,\mu}$-module $M$ we consider the subspace
\begin{align}\label{eqn:M[j]}
	M[j]\coloneqq M_{|}[j].
\end{align}
This somewhat redundant notation is justified by the following fact.
\begin{proposition}\label{pro:j-unico}
		Let $M$ be an $\A_{\lambda,\mu}$-module. For each $j\in\J_0$, $M[j]\subseteq M$ is an $\mA_{\lambda,\mu}$-submodule and 
		$M$ decomposes as $M=\bigoplus_{j=0}^{\ell-1}M[j]$. 	
	\end{proposition}
	Hence, if $M$ is indecomposable, then there is $j$ so that $M=M[j]$.
\pf
Let $x\in M[j]$ and write $y=a_0\cdot x$ as $y=\sum_{i=0}^{\ell-1} y_i$, $y_i\in M_{|}[i]$. As $t^2a_0=a_0t^2$ in $\A_{\lambda,\mu}$, the result follows by Lemma \ref{lem:uniquej}. Indeed, $t^2a_0\cdot x=t^2\cdot y=\sum_{i=0}^{\ell-1}\zeta^i\, y_i$ and $a_0t^2\cdot x=\zeta^j y=\sum_{i=0}^{\ell-1}\zeta^j\, y_i$.
\epf 	
We can decompose the regular representation as $\mA_{\lambda,\mu}\simeq \bigoplus_{j=0}^{\ell-1} \A[j]$, where $\A[j]\coloneqq \A(S_j^+)\oplus \A(S_j^-)\oplus \A(M_j)$.	We thus get that $\dim \A[j]=72$ for each $j\in\J_0$ and a $j$-form of formula \eqref{eqn:dim-simple-proy}:
	\begin{align}\label{eqn:dim-simple-proy-j}
		72=\sum_{L=L[j]}\dim L\dim P(L).
	\end{align}

\subsubsection{Sums of one-dimensional modules}\label{subsec:powersof1}

We look at modules $M$ so that $M_{|}\simeq S_{j}^{+p}\oplus S_j^{-q}$.  
By Proposition \ref{pro:j-unico}, $a_0\cdot S_{j}^{+p}\subset S_{j}^{-q}$, and vice versa. 

\begin{remark}\label{rem:suma-dim1}
	Let $j\in\J$, $p,q>0$ and assume there is an $\A_{\lambda,\mu}$-module $M$ so that $M_{|}\simeq S_{j}^{+p}\oplus S_j^{-q}$. 
Then $	\mu=\frac\lambda3\neq0$.

Indeed, we have that the action of $s$ is trivial: $[s]=\Id$ and therefore $[a_0]=[a_1]=[a_2]$.
Hence, from \eqref{eqn:rels_def} we get
$a_0^2=\mu(1-t^2)$ and $3a_0^2=\lambda(1-t^2)$.	As $[t]^2\neq \Id$, we obtain the equality.\qed
\end{remark}

This allows us to introduce a simple module in this setting.

\begin{lemma}\label{lem:tj}
Assume $\mu=\frac{\lambda}{3}$. Then for each $j\in\J$ there is a simple $\mA$-module $T_j$ so that $T_{j|}\simeq S_j^+\oplus S_j^-$ and the action of $a_0$ is determined by 
\begin{align}\label{eqn:Tj}
[a_0]=\begin{pmatrix*} 0&\mu(1-\zeta^j)\\1&0\end{pmatrix*}.
\end{align}
\end{lemma}

For the case $j=0$ this shifts from irreducibility to indecomposability.

\begin{lemma}\label{lem:t0}
There are, up to isomorphism, two indecomposable non-simple $\mA$-modules $T_0^+$ and $T_0^-$ with $(T_{0}^\pm)_{|}\simeq S_0^+\oplus S_0^-$. The action of $a_0$ on $T_0^+$ and $T_0^-$ is determined by 
\begin{align}\label{eqn:t0}
[a_0]&=\begin{pmatrix*} 0&0\\1&0\end{pmatrix*},  & [a_0]&=\begin{pmatrix*} 0&1\\0&0\end{pmatrix*}
\end{align}	
respectively. In other words, $L_0^\mp$ is a submodule of $T_0^\pm$, which is generated by the component $L_0^\pm$; that is, this module fits in the extension
\[
0\longrightarrow L_0^\mp  \longhookrightarrow T_0^\pm \longtwoheadrightarrow L_0^\pm  \longrightarrow 0.
\]
\end{lemma}
\pf
Straightforward.
\epf
We shall make use of the following immediate fact along the proof of Theorem \ref{thm:simple} in \S\ref{sec:classification}.
\begin{remark}\label{rem:dim1-notprojective}
An easy consequence of Lemma \ref{lem:t0} is the realization that the one-dimensional modules $L_0^\pm$  are not projective, and in particular their projective covers $P(L_0^\pm)$ satisfy $\dim P(L_0^\pm)\geq 2$. 
\end{remark}

The modules form Lemmas \ref{lem:tj} and \ref{lem:t0} characterize simple and indecomposable modules of such shape.

\begin{lemma}\label{lem:suma-dos-dim1}
Let $M$ be an $\mA$-module with $M_{|}\simeq S_{j}^{+p}\oplus S_j^{-q}$, $j\in\J_0$.
	\begin{enumerate}[leftmargin=*]
\item If $j\in\J$, then $\mu=\frac{\lambda}{3}$ and $p=q$. In particular, $M\simeq (T_j)^p$ and $M$ is indecomposable if and only if it is simple. 
\item If $j=0$, then $M\simeq L_0^{+a}\oplus L_0^{-b}\oplus T_0^{+c}\oplus T_0^{-d}$, with $a+c+d=p$, $b+c+d=q$. If $M$ is indecomposable but not simple, then $M\simeq T_0^\pm$.
	\end{enumerate}
\end{lemma}
\pf
Assume $j\in\J$ and $p\leq q$. Set $x=x_1$ in $M$ so that $\lg x_1\rg_{|}\simeq S_j^+$. By Lemma \ref{lem:nontrivial action dim1} we have
$0\neq y_1\coloneqq a_0\cdot x_1$. Moreover, $\lg y_1\rg_{|}\simeq S_j^-$. This also gives
\[
a_0\cdot y_1=a_0^2\cdot x_1=\mu(1-\zeta^j)x_1.
\] 
Which shows that $T_j\simeq \lg x_1,y_1\rg\subset M$. We can repeat this process for each element in a basis $\{x_1,x_2,\dots, x_p\}$ in the component $S_j^{+p}$. We claim the the corresponding set $\{y_1,y_2,\dots, y_p\}$ is linearly independent. Indeed, if $\sum_{i=0}^p c_i y_i=0$, then $0=\sum_{i=0}^p c_i a_0\cdot y_i=\mu(1-\zeta^j)\sum_{i=0}^p c_i x_i$, which is a contradiction. Now, if $q>p$ and $y\notin \lg y_1,y_2,\dots, y_p\rg$ then $0\neq x=a_0\cdot y$ and $\lg x\rg_{|}\simeq S_j^+$. But this gives $\mu(1-\zeta^j)y=a_0\cdot x\in\lg y_1,y_2,\dots, y_p\rg$, again a contradiction. Thus $p=q$ and $M\simeq T_j^p$. As $M$ is indecomposable, $p=1$ and \textit{(1)} follows.

Let $\{x_1,\dots,x_p\}$ and $\{y_1,\dots,y_q\}$  bases of the components $S_0^+$ and $S_0^-$ respectively. If $a_0\cdot x_i=0$ for each $i=1,\dots,p$ and $a_0\cdot y_i=0$, for all $i=1,\dots q$, then $M\simeq L_0^{+p}\oplus L_0^{-q}$. Otherwise, assume that there are $p'\leq p$, $q'\leq q$, not simultaneously zero, so that $a_0\cdot x_i\neq 0$ for $i\in X=\{1,\dots,p'\}$ and $a_0\cdot y_i\neq 0$, for $i\in Y=\{1,\dots, q'\}$. 

It follows that $y_i\notin \lg a_0\cdot x_1,\dots,a_0\cdot x_{p'}\rg$, $i\in Y$ and  $x_i\notin \lg a_0\cdot y_1,\dots,a_0\cdot y_{q'}\rg$, $i\in X$. Indeed, $a_0\cdot \lg a_0\cdot x_1,\dots,a_0\cdot x_{p'}\rg=0$ and $a_0\cdot \lg a_0\cdot y_1,\dots,a_0\cdot y_{q'}\rg=0$. 

Hence $p'+q'\leq p$ and $S_0^{+p}\supseteq \lg x_i,a_0\cdot y_j|i\in X,j\in Y \rg$ is a linearly independent set. Similarly,  $p'+q'\leq q$ and $\lg y_i,a_0\cdot x_j|i\in Y,j\in X \rg \subset S_0^{-q}$.  

As well, $\lg x_i,a_0\cdot x_i\rg\simeq T_0^+\subset M$, $i\in X$, and  $\lg y_i,a_0\cdot y_i\rg\simeq T_0^-\subset M$, $i\in Y$. Namely $T_0^{+p'}\oplus T_0^{-q'}\subset M$.

If $p'+q'<p$, then there are elements $x_k$, $k=p'+q'+1,\dots, p$ so that $a_0\cdot x_k=0$ and $x_k\notin a_0\cdot S_0^-$. Similarly, if $p'+q'<q$, then  there are elements $y_k\notin a_0\cdot S_0^+$ so that $a_0\cdot y_k=0$, $k=p'+q'+1,\dots,q$ and $S_0^{-q}=\lg  y_1,\dots,y_{q'},a_0\cdot x_1,\dots,a_0\cdot x_{p'},y_{p'+q'+1},\dots y_q\rg$.

Setting $c=p'$, $d=q'$, $a=p-p'+1$, $b=q-q'+1$ we obtain \textit{(2)}.
\epf

\subsubsection{Powers of the two-dimensional module}\label{subsec:powersof2}

Now we study modules sustained on sums $M_j\oplus\dots\oplus M_j$. We show that they are, generically, sums of simple modules $L_j(a,c)$ as in Lemma \ref{lem:LMj}.

\begin{proposition}\label{pro:sumadedim2}
Assume there is a $\A_{\lambda,\mu}$-module $M$ such that $M_{|}\simeq M_j^r$. Then
$\mu=\frac{\lambda}{3}$ or $j=0$. 
\begin{enumerate}[leftmargin=*]
\item $M$ is simple if and only if $r=1$ and $M\simeq L_j(a,c)$, $a,c\in\{\pm\sqrt{\lambda/3}\}$.
\item If $\lambda\neq 0$, then $M$ is a direct sum of simple modules of dimension 2. 
\end{enumerate}
\end{proposition}
\pf
Set $\{v_1, w_1\}, \{v_2, w_2\}, \dots, \{v_r, w_r\}$  bases for each copy of $M_j$. We fix the basis $\{v_1, v_2,\dots, v_r, w_1, w_2, \dots, w_r\}$ for $M$. Let  $v=(v_1, \dots, v_r)^t$, $w=(w_1, \dots, w_r)^t$  and $\alpha=(\alpha_{i,j})$, $\beta=(\beta_{i,j}))\in\k^{l\times l}$ be such that 
$
a_0\cdot v= \alpha v + \beta w$.
That is $a_0\cdot v_i=\sum \alpha_{i,j}v_j+\sum \beta_{i,j}v_j$. Thus $
	a_0\cdot w=a_0\cdot t	\cdot v= -\beta\zeta^j v - \alpha w$.
	
By \eqref{eqn:rels_def}, $a_0^2\cdot v = \mu(1-\zeta^j)v$ while 
\begin{align*}
	a_0^2\cdot v &=\alpha a_0 v + \beta a_0 w=\alpha (\alpha v + \beta w) + \beta (-\beta\zeta^j v - \alpha w)\\
	&=(\alpha^2-\beta^2\zeta^j) v+ (\alpha\beta-\beta\alpha)w.
\end{align*}
Therefore, we obtain the relations
\begin{align}\label{eqn:alphabeta}
	\alpha\beta=\beta\alpha, && \alpha^2-\beta^2\zeta^j=\mu(1-\zeta^j)\id.
\end{align}
On the other hand, as $a_1=sa_0s^{-1}$, $a_2=s^{-1}a_0s$,
\begin{align*}
		a_1\cdot v&=\alpha v + \xi\beta w, & 		a_1 \cdot w&=-\beta\zeta^j\xi^2 v-\alpha w,\\
		a_2 \cdot v&=\alpha v + \xi^2\beta w, & 		a_2 \cdot w&=-\beta\zeta^j\xi v-\alpha w.
\end{align*}
Thus $(a_0a_1+a_1a_2+a_2a_0)\cdot v=3(\alpha^2-\xi\zeta^j\beta^2) v$.
While $\lambda(1-st^2)\cdot v=\lambda(1-\xi\zeta^j)v$. On the other hand, $\lambda(1-st^2)\cdot w= \lambda(1-\xi^2\zeta^j)w$ and 
	$(a_0a_1+a_1a_2+a_2a_0)\cdot w=	3(\alpha^2-\xi^2\zeta^j\beta^2) w$.
Thus
\begin{align*}
		3(\alpha^2-\xi\zeta^j\beta^2)&=\lambda(1-\xi\zeta^j)\id, &
		3(\alpha^2-\xi^2\zeta^j\beta^2)&=\lambda(1-\xi^2\zeta^j)\id.
	\end{align*}
Hence $\alpha^2=\beta^2=\frac{\lambda}{3}\id$. If we combine this with the second identity in \eqref{eqn:alphabeta} we get that either $\mu=\frac{\lambda}{3}$ or $j=0$, as stated.

If $\lambda\neq0$, then the identities $\alpha^2=\beta^2=\frac{\lambda}{3}\id$ together with the fact that $\alpha\beta=\beta\alpha$ by \eqref{eqn:alphabeta}, imply that $\alpha$ and $\beta$ are simultaneously diagonalizable and thus there is a basis in which $M$ is a direct sum of modules of dimension 2. Indeed, there is an invertible matrix $p$ and diagonal matrices $d=\diag(d_1,\dots,d_r)$ and $d'=\diag(d'_1,\dots,d'_r)$ such that \begin{align*}
	p\alpha p^{-1}=d,&& p\beta p^{-1}=d'.
\end{align*}
Hence, if $v_i'=\sum p_{ij}v_j$, $w_i'=\sum p_{ij}w_j$, 
then $a_0\cdot v_i'=d_iv_i'$ y $a_0\cdot w_i'=d'_iw_i'$. 

Furthermore, note that $t\cdot v_i'=w_i'$ and $t\cdot w_i'=\zeta^jv_i'$. Analogously, $s\cdot v_i'=\xi v_i'$ and $s\cdot w_i'=\xi^2 w_i'$. 
Thus, for each $i$, $\k\{v_i',w_i'\}$ is a 2-dimensional submodule of $M$ and $M$ is completely reducible as a sum of modules of dimension 2 and therefore it it not simple if $r>1$.

When $\lambda=0$, then we necessarily get $j=0$ as $\mu\neq 0$, from \eqref{eqn:alphabeta}, since the commuting matrices $\alpha$ y $\beta$ satisfy $\alpha^2=\beta^2=0$. As a result, they are simultaneously triangularizable. 
In particular, it implies the existence of $v_1'$ and $w_1'=t\cdot v_1'$ such that $a_0\cdot v_1'=a_0\cdot w_1'=0$, which defines a submodule $\k\{v_1', w_1'\}\subset M$. 
Thus $M$ is not simple if $r>1$.
\epf

\begin{example}\label{exa:extensionM0}
Fix $j=0$. Assume $\lambda=0$ and let $\lg v,w\rg$ and $\lg v',w'\rg$ be two copies of the $G$-module $M_0$. Then 
\begin{align*}
a_0\cdot v&=v', & a_0\cdot v'&=0, & a_0\cdot w&=w', & a_0\cdot w'&=0
\end{align*}
define  a non-simple, indecomposable module $M$ with $M_{|}\simeq M_0^2$. In particular, such $M$ defines an extension
\[
0\to M_0\hookrightarrow M\twoheadrightarrow M_0\to 0.
\]
See \S\ref{sec:extensions} for a complete description of extension between simple $\mA$-modules.
\end{example}

\subsubsection{Very mixed modules}\label{subsec:verymixed}

Let $L$ be an $\mA_{\lambda,\mu}$-module so that $L=L[j]$ for some $j\in\J_0$. Recall that this is the case when $L$ is indecomposable.
In particular, this fixes  integers $p,q,r\geq 0$ so that
\begin{align}\label{eqn:L-decomp}
	L_{|}\simeq (S^{+}_{j})^{p}\oplus (S^{-}_{j})^{q}\oplus M^r_{j}.
\end{align}
Associated to this we introduce the following scalars:
	\begin{align}
		\begin{split}\label{eqn:scalars-mixed}
	f_j=\dfrac{1}{9b}(\lambda-3\mu)(1-\zeta^j),  \qquad  g_j=\frac{1}{9}(3\mu+2\lambda)(1-\zeta^j),\\
f_{\pm j}=\dfrac{\mu}{3}(1-\zeta^{\pm j})+\dfrac{\lambda}{9}(2+\zeta^{\pm j}).\qquad \qquad 
		\end{split}
\end{align}
We omit the subscript $j$ when it is clear, we write $b=b_j$, $g=g_j$, etc.
\begin{lemma}\label{lem:general-a0}
	Let $L$ be an $\A_{\lambda,\mu}$-module as in \eqref{eqn:L-decomp}.
	Then there matrices $\alpha\in\k^{q\times p}$, $\beta\in\k^{p\times q}$, $\delta_{+}\in\k^{p\times r}$, $\delta_{-}\in\k^{q\times r}$, $\eta\in\k^{r\times p}$, $\sigma\in\k^{r\times q}$ and $\theta,\tau\in\k^{r\times r}$ subject to conditions:
	\begin{align}\label{eqn:general1}
		\begin{split}
		\beta\alpha&=g\,\id_p, \qquad  \alpha\beta=g\,\id_q,\\
		\theta^2&=f_+\, \id_r  \qquad  \tau^2=f_-\,\id_r.
		\end{split}
	\end{align}
	\begin{align}\label{eqn:general2}
		\begin{split}
		&	\begin{cases}
			\delta_{+}\sigma=\delta_{-}\eta=0,\\
			\delta_+\eta=-\delta_{-}\sigma=f\,\id_p,
		\end{cases} \qquad \begin{cases}
			\eta\delta_{+}-\sigma\delta_{-}=f\,\id_r,\\
			\eta\delta_{+} +\sigma\delta_{-}= \tau\theta-\theta\tau,
		\end{cases}
		\end{split}
		\end{align}
	\begin{align}\label{eqn:general3}
	\begin{split}
		\begin{cases}
			\beta\delta_{-}-  \delta_+(b\tau +\theta)=0=\alpha\delta_{+} + \delta_{-}(b\tau- \theta),\\
			\sigma\alpha - (b\tau-\theta)\eta=0=\eta\beta + (b\tau+\theta)\sigma,
		\end{cases}
		\end{split}
\end{align}
	so that the action of $a_0$ is described by the matrix (notice the transpose):
	\begin{align}\label{eqn:a0matrix}
		[a_0]=\begin{pmatrix}
			0&  \beta&  -b  \delta_+&\delta_+\\
			  \alpha&0& b  \delta_- &\delta_-\\
			\eta&\sigma&  \theta&\tau\\
			  -b   \eta&   b  \sigma&-\zeta^j  \tau &-  \theta
		\end{pmatrix}^t.
	\end{align}

Conversely, any collection of matrices $(\alpha,\beta,\gamma,\delta_{\pm},\eta,\sigma,\theta,\tau)$ so that \eqref{eqn:general1}, \eqref{eqn:general2} and \eqref{eqn:general3} hold determines an $\mA_{\lambda,\mu}$-module $L_j^{p,q,r}(\alpha,\beta,\delta_{+},\delta_-,\eta,\sigma,\theta,\tau)$ so that \eqref{eqn:L-decomp} holds and the action of $a_0$ is via \eqref{eqn:a0matrix}.
\end{lemma}

\pf
Let us fix generators 
\begin{align*}
x&=(x_1, \dots, x_{p})^t, & y&=(y_1, \dots, y_{q})^t, & \{v=(v_1, \dots, v_{r})^t, w=(w_1, \dots, w_{r})^t\}
\end{align*} 
of $(S^{+}_{j})^{p}, (S^{-}_{j})^{q}$ y $M^r_{j}$, respectively. By abuse of notation, we keep the notation $b$ for the scalar so that $b^2=\zeta^j$ as well as for the matrices $b\coloneqq\diag(b)$ both in $\k^{p\times p}$ and $\k^{q\times q}$. 

We fix matrices $\alpha_{+}\in\k^{p\times p}$, $\alpha_{-}\in\k^{q\times q}$, $\beta_{+}\in\k^{p\times q}$, $\beta_{-}\in\k^{q\times p}$,  together with $\gamma_{+},\delta_+\in\k^{p\times r}$ and $\gamma_{-}, \delta_{-}\in\k^{q\times r}$ so that
\[
a_0\cdot x= \alpha_{+} x + \beta_{+} y + \gamma_{+} v + \delta_{+} w, \qquad a_0\cdot y= \alpha_{-} x + \beta_{-} y+ \gamma_{-} v + \delta_{-} w,
\]
First, as $a_0t+ta_0=0$ we get $\alpha_{+}=0=\beta_{-}$. Let us rename $\alpha_{-}=\alpha$ y $\beta_{+}=\beta$. Then,
\[
a_0\cdot x= \beta y + \gamma_+ v + \delta_+ w, \qquad a_0\cdot y= \alpha x + \gamma_- v + \delta_- w.
\]
From this identity we also arrive to
\begin{align*}
	\left\{ \begin{array}{lc}
		a_0\cdot t\cdot x=b a_0\cdot x= b\beta y + b\gamma_+ v + b\delta_+ w\\
		-t\cdot a_0\cdot x=b\beta y - \gamma_+ w - \zeta^j\delta_+ v
	\end{array}
	\right.
	\implies b\delta_+=-\gamma_+.
\end{align*}
and analogously $\gamma_{-}=b\delta_{-}$. Thus we end up with
\begin{align*}
	a_0\cdot x= \beta y + \delta_+ (w-bv), \qquad a_0\cdot y= \alpha x + \delta_- (w+bv).
\end{align*}
On the other hand, let $\eta\in\k^{r\times p}$, $\sigma\in\k^{r\times q}$, and $\theta,\tau\in\k^{r\times r}$ be such that
\begin{align*}
	a_0\cdot v&=\eta x + \sigma y +\theta v + \tau w.
\end{align*}
It follows that
\begin{align*}
	a_0\cdot w&=a_0\cdot t\cdot v=-t\cdot a_0\cdot v= -b\eta x + b\sigma y - \zeta^j\tau v -\theta w.
\end{align*}
Hence $[a_0]$ is as in  \eqref{eqn:a0matrix}.
From $a_1=sa_0s^2$ and $a_2=s^2a_0s$:
\begin{align*}
[a_1]&=\begin{psmallmatrix}
	0&  \beta& -b\xi  \delta_+&\xi^2  \delta_+\\
	  \alpha&0&b\xi  \delta_-&\xi^2  \delta_-\\
	\xi^2  \eta&\xi^2  \sigma&  \theta&\xi   \tau\\
	-b\xi   \eta&b\xi  \sigma&-\zeta^j\xi^2  \tau&-  \theta
\end{psmallmatrix}^t, &&
[a_2]&=\begin{psmallmatrix}
	0&  \beta&-b\xi^2  \delta_+&\xi  \delta_+\\
	  \alpha&0&b\xi^2  \delta_-&\xi \delta_-\\
	\xi  \eta&\xi  \sigma&  \theta&\xi^2   \tau\\
	-b\xi^2   \eta&b\xi^2  \sigma&-\zeta^j\xi  \tau&-  \theta
\end{psmallmatrix}^t.
\end{align*}
Namely,
\begin{align*}
	a_1\cdot x &=\beta y-b\xi\delta_+v+\xi^2\delta_+w, && a_2\cdot x=\beta y-b\xi^2\delta_+v+\xi\delta_+w\\
	a_1\cdot y &= \alpha x+ b\xi\delta_-v+\xi^2\delta_-w,&& a_2\cdot y=\alpha x+ b\xi^2\delta_-v+\xi\delta_-w,\\
	a_1\cdot v &= \xi^2\eta x + \xi^2\sigma y + \theta v + \xi\tau w,&& a_2\cdot v=\xi\eta x + \xi\sigma y + \theta v + \xi^2\tau w.
\end{align*}
Let us analyze the identity $[a_0]^2=\mu(\id-[t]^2)=\mu(1-\zeta^j)\id$:
\begin{align*}
	a_0^2\cdot x&=(\beta\alpha - 2 b\delta_+ \eta)x +(b\beta\delta_- -b\delta_+\theta - \zeta^j\delta_+\tau)v\\
	&+ (\beta\delta_{-}-  b\delta_+\tau -\delta_+\theta)w=\mu(1-\zeta^j)x,\\
	a_0^2\cdot y&=(\alpha\beta + 2 b\delta_{-}\sigma)y + (-b\alpha\delta_{+} + b\delta_{-}\theta - \zeta^j\delta_{-}\tau)v\\
	&+ (\alpha\delta_{+} + b\delta_{-}\tau - \delta_{-}\theta)w=\mu(1-\zeta^j)y,\\
	a_0^2\cdot v&=(\sigma\alpha + \theta\eta -b\tau\eta )x + (\eta\beta + \theta\sigma + b\tau\sigma)y\\
	&+(-b\eta\delta_{+} + b\sigma\delta_{-} + \theta^2 -\zeta^j\tau^2)v + (\eta\delta_{+} + \sigma\delta_{-} + \theta\tau -\tau\theta )w\\
	&=\mu(1-\zeta^j)v,\\
	a_0^2\cdot w&=(b\sigma \alpha -\zeta^j\tau \eta + b\theta \eta)x + (-b\eta \beta -\zeta^j\tau\sigma -b\theta\sigma)y\\
	&+\zeta^j(\eta\delta_{+} +\sigma\delta_{-} -\tau\theta +\theta\tau)v + (-b\eta\delta +b\sigma\delta_{-}-\zeta^j\tau^2 + \theta^2)w\\
	&=\mu(1-\zeta^j)w.
\end{align*}
Thus we obtain $\beta \alpha - 2 b\delta_+ \eta=\alpha\beta + 2 b\delta_{-}\sigma=\mu(1-\zeta^j)$ and  $b(\sigma\delta_{-}-\eta\delta_{+}) + \theta^2 -\zeta^j\tau^2=\mu(1-\zeta^j)$, together with
\begin{align*}
	\beta\delta_{-}-  b\delta_+\tau -\delta_+\theta&=0, & \alpha\delta_{+} + b\delta_{-}\tau- \delta_{-}\theta&=0,\\
	\sigma\alpha + \theta\eta -b\tau\eta&=0, & \eta\beta + \theta\sigma + b\tau\sigma&=0,\\
	\eta\delta_{+} + \sigma\delta_{-} + \theta\tau -\tau\theta&=0.
\end{align*}
On the other hand, we have that
\[(a_0a_1+a_1a_2+a_2a_0)\cdot(x, y, v, w)=\lambda\big((1-\zeta^j)x,(1-\zeta^j)y,(1-\xi\zeta^j)v,(1-\xi^2\zeta^j)w\big)\]
which is 
\begin{align*}
&3(\beta\alpha +b\delta_{+}\eta)x + 3(\xi^2-\xi)\delta_{+}\sigma y=\lambda(1-\zeta^j)x,\\
&b(\xi-\xi^2)\delta_{-}\eta x+3(\alpha\beta -b\delta_{-}\sigma)y=\lambda(1-\zeta^j)y\\
&3(\theta^2 -\zeta^j\xi\tau^2 +b\xi^2(\sigma\delta_{-}-\eta\delta_{+})v=\lambda(1-\xi\zeta^j)v\\
&3(\theta^2-\zeta^j\xi^2\tau^2+b\xi(\sigma\delta_{-}-\eta\delta_{+}))w=\lambda(1-\xi^2\zeta^j)w
\end{align*}
and therefore $\delta_{+}\sigma=0$ and $\delta_{-}\eta=0$. As well,
\begin{align*}
	\beta\alpha +b\delta_{+}\eta=\alpha\beta -b\delta_{-}\sigma&=\frac{1}{3}\lambda(1-\zeta^j),\\
	\theta^2 -\zeta^j\xi\tau^2 +b\xi^2(\sigma\delta_{-}-\eta\delta_{+})&=\frac{1}{3}\lambda(1-\xi\zeta^j),\\
	\theta^2-\zeta^j\xi^2\tau^2+b\xi(\sigma\delta_{-}-\eta\delta_{+})&=\frac{1}{3}\lambda(1-\xi^2\zeta^j).
\end{align*}
Hence we obtain the systems:
\begin{align}\label{eqn:sistema-alphabeta}
	&\begin{cases}
		\beta\alpha +b\delta_{+}\eta=\frac{1}{3}\lambda(1-\zeta^j),\\
		\beta \alpha - 2 b\delta_+ \eta=\mu(1-\zeta^j),
	\end{cases} &&
	&\begin{cases}
		\alpha\beta -b\delta_{-}\sigma=\frac{1}{3}\lambda(1-\zeta^j),\\
		\alpha\beta + 2 b\delta_{-}\sigma=\mu(1-\zeta^j),\\
	\end{cases} 
\end{align}
together with
\begin{align*}
	\begin{cases}
		\theta^2 -\zeta^j\tau^2+b(\sigma\delta_{-}-\eta\delta_{+}) =\mu(1-\zeta^j)\\
		\theta^2 -\zeta^j\xi\tau^2 +b\xi^2(\sigma\delta_{-}-\eta\delta_{+})=\frac{1}{3}\lambda(1-\xi\zeta^j),\\
		\theta^2-\zeta^j\xi^2\tau^2+b\xi(\sigma\delta_{-}-\eta\delta_{+})=\frac{1}{3}\lambda(1-\xi^2\zeta^j).
	\end{cases} 
\end{align*}
A straightforward resolution of these systems gives the conditions in the statement of the lemma.

Conversely, if $[a_0]$ is as in \eqref{eqn:a0matrix} and equations \eqref{eqn:general1}--\eqref{eqn:general3} hold, then it is straightforward to check that relations
\eqref{eqn:rels_def} hold, for $[a_1]=[s][a_0][s]^2$ and $[a_2]=[s]^2[a_0][s]$ as above.
\epf

\begin{corollary}\label{cor:dim4}
	Fix $j\in\J_0$ and assume there is a $\A_{\lambda,\mu}$-module  $M$ so that $M_{|}\simeq S^{+}_{j}\oplus S^{-}_{j}\oplus M_{j}$. Then $M$ is not simple.
	
	In particular, there are no simple modules of dimension 4.
\end{corollary}

\pf
We fix, as above, a basis $\{x,y,v,w\}$ for $M$, in such a way that $\{x\}$, $\{y\}$ and $\{v,w\}$ are the canonical generators of $S^{+}_{j}, S^{-}_{j}$ y $M_{j}$, respectively. 
Set $b\coloneqq b_j^+$. We have that
\begin{align*}
	[a_0]=\begin{psmallmatrix}
		0&\alpha&\eta&-b \eta\\
		\beta&0&\sigma&b\sigma\\
		b\delta_+&b\delta_-&\theta&-\zeta^j\tau\\
		\delta_+&\delta_-&\tau&-\theta
	\end{psmallmatrix}.
\end{align*}
We invoke  \eqref{eqn:general2} where $\delta_{+}\sigma=\delta_{-}\eta=0.$ and analyze the four possibilities:
\begin{align*}
	(i)\,\delta_-=\delta_+&=0, & 
	(ii)\,\eta=\delta_+&=0, & 
	(iii)\,\delta_-=\sigma&=0, & 
	(iv)\,\eta=\sigma&=0.
\end{align*}
Now, in (i) the span of $\{x,y\}$ determines a submodule $\simeq S_j^+\oplus S_j^-$ while in (iv) the span $\{v,w\}$ determines a submodule $\simeq M_j$.
For (ii), we use 
\[
\beta\delta_--\delta_+(b\tau+\theta)=0\Rightarrow \beta\delta_-=0.
\]
As we can assume $\delta_-\neq 0$ (otherwise this is case (i)) we get $\beta=0$ and thus $\lg x\rg$ defines a submodule $\simeq S_j^+$. 
Finally, in case (iii) we use:
\[
\alpha\delta_+ +\delta_-(b\tau -\theta)=0\Rightarrow \alpha\delta_+=0.
\]
Assuming $\delta_+\neq 0$, we get $\alpha=0$ and  thus $\lg y\rg$ defines a submodule $\simeq S_j^-$. 
\epf

\begin{lemma}\label{lem:1+2-notsimple}
Let $M$ be a $\mA_{\lambda,\mu}$-module so that either 
\[
M_{|}\simeq (S_j^+)^p\oplus M_j^r \quad \text{ or } \quad M_{|}\simeq (S_j^-)^q\oplus M_j^r
\]
with $p,q,r\geq 1$. Then $M$ is not simple.

In particular, there are no simple modules of dimension 3 or 4.
\end{lemma}

\pf
Assume $M_{|}\simeq (S_j^+)^p\oplus M_j^r$, the other case is analogous. We have that
$	[a_0]=\begin{psmallmatrix}
		0&-b  \delta_+&\delta_+\\
		\eta&  \theta&\tau\\
		  -b   \eta&  -\zeta^j  \tau&-  \theta
	\end{psmallmatrix}^t$, 
with $\delta_+(b\tau+\theta)=0$, $\tau\theta-\theta\tau=f\,\id$ and 
\begin{align*}
\theta^2&=f_+\,\id, & \tau^2&=f_-\,\id, & \delta_+\eta&=f\,\id, & \eta\delta_+&=f\,\id.
\end{align*}
First, assume $\lambda\neq 3\mu$, so $f\neq 0$. In particular, $\delta_+$ is invertible (and $p=r$) and thus $b\tau+\theta=0$.
But this implies $\tau\theta-\theta\tau=0$, a contradiction.

Now, if $\lambda= 3\mu$, then $\delta_+\eta=0=\eta\delta_+$. Hence ${}^t\delta_+$ is not injective and we can rearrange the basis of $M_j^r$ to get $a_0\cdot x_1=0$ and thus $\lg x_1\rg\subset M$.

The last assertion follows since there are no simple modules supported over sums of one-dimensional $\G_{3,\ell}$-modules nor over $M_j^2$, and we have just discarded $S_j^{\pm2}\oplus M_j$.
\epf

\begin{lemma}\label{lem:mixed-Mj-1}
	Let $M$ be a $\mA_{\lambda,\mu}$-module so that
	\[
	M_{|}\simeq (S_j^+)^p\oplus  (S_j^-)^q\oplus M_j
	\]
	with $p,q\geq 1$. Then $M$ is not simple. In particular, there are no simple modules of dimension 5.
\end{lemma}
\pf
In this case we get
\begin{align*}
\eta\delta_+&=\frac{1}{2}f, & 
\sigma\delta_-&=-\frac{1}{2}f, & 
\delta_+\eta&=\frac{1}{2}f, & 
\delta_-\sigma&=-f.
\end{align*}
If $f\neq 0$ this implies $p=q=1$ and Corollary \ref{cor:dim4} applies. If $f=0$, then we can assume, by rearranging the basis, that $a_0\cdot x_1\in\lg y_1\rg$
and therefore $\lg x_1,y_1\rg$ determines a submodule.
\epf

\subsubsection{A six-dimensional simple module for $j\in\J$}\label{subsec:verymixed6}

We define a 6-dimensional simple module $N_j$ for each $j\in\J$ when $3\mu\neq\lambda$. We show that this is not possible when $3\mu = \lambda$.
Recall the scalars from \eqref{eqn:scalars-mixed} and let us add: 
\begin{align*}
	h_j&\coloneqq f_{+j}+b_j^2f_{-j}+g_j=\frac\lambda9(5+b_j^2)+\frac\mu3(1-b_j^2).		
\end{align*}
As well, we shall consider the following set.
\begin{definition}\label{def:solutions}
For each $j\in\J$ we let $\Up_j\subset \k^4$ be the subset of 4-tuples $\up=(c_1,c_2,t_1,t_2)$ such that
\begin{align}\label{eqn:conditions-six}
	c_1c_2&=g_j, & t_1t_2&=f_{+j}, & c_1t_1+c_2t_2=h_j.
\end{align}
\end{definition}

It is easy to see that $\Up_j\neq\emptyset$ for each $j\in\J$. We shall give a complete description of these sets in  \S\ref{sec:isoclasses}, in order to deal with the isomorphism classes of the 6-dimensional simple modules to be defined in the proposition that follows.

\begin{proposition}\label{pro:dimension6}
	Assume $3\mu\neq\lambda$, $j\in\J$ and let $\up=(c_1,c_2,t_1,t_2)\in\Up_j$. 
	There is an irreducible $\A_{\lambda, \mu}$-module $N_j(\up)=\lg x,y,v_1,v_2,w_1,w_2\rg$ so that
	\begin{align}\label{eqn:isotypic}
		\lg x \rg&\simeq S_j^+, & 
		\lg y \rg&\simeq S_j^-, &
		\lg v_1,w_1 \rg&\simeq \lg v_2,w_2 \rg\simeq M_j.
	\end{align}
	and the action of $a_0$ on this basis is codified by the matrix	\begin{align}\label{eqn:dimension6}
		[a_0]=\begin{psmallmatrix}
			0&c_2&f&0&-bf&0\\
			c_1&0&0&-f&0&-bf\\
			-b&0&0&t_1&0&b(c_2-t_1)\\
			0&b&t_2&0&b(t_2-c_1)&0\\
			1&0&0&(t_1-c_2)/b&0&-t_1\\
			0&1&(c_1-t_2)/b&0&-t_2&0
		\end{psmallmatrix}.
	\end{align}
	Conversely, if $N$ is an irreducible 6-dimensional $\mA_{\lambda,\mu}$-module, then there are $j\in\J$ and $\up\in\Up_j$ so that $N\simeq N_j(\up)$.
\end{proposition}

\pf
Set $N=N_j$ and let $s,t\in\k^{6\times 6}$ denote wlog the matrices corresponding to elements $s,t\in\G_{3,\ell}$, namely:
\begin{align*}
	[t]=\begin{psmallmatrix}
		b&0&0&0&0&0\\
		0&-b&0&0&0&0\\
		0&0&0&0&b^2&0\\	
		0&0&0&0&0&b^2\\
		0&0&1&0&0&0\\
		0&0&0&1&0&0
	\end{psmallmatrix},
	&&
	[s]=\diag(1,1,\xi,\xi,\xi^2,\xi^2).
\end{align*}
It easily follows that $-[t][a_0]=[a_0][t]$. From $[a_1]=s[a_0]s^2$ y $[a_2]=s^2[a_0]s$:
\begin{align*}
	[a_1]&=\begin{psmallmatrix}
		0&c_2&\xi^2f&0&-b\xi f&0\\
		c_1&0&0&-\xi^2f&0&-b\xi f\\
		-b\xi&0&0&t_1&0&b\xi^2(c_2-t_1)\\
		0&b\xi&t_2&0&b\xi^2(t_2-c_1)&0\\
		\xi^2&0&0&\xi(t_1-c_2)/b&0&-t_1\\
		0&\xi^2&\xi(c_1-t_2)/b&0&-t_2&0
	\end{psmallmatrix},\\
	[a_2]&=\begin{psmallmatrix}
		0&c_2&\xi f&0&-b\xi^2 f&0\\
		c_1&0&0&-\xi f&0&-b\xi^2 f\\
		-b\xi^2&0&0&t_1&0&b\xi(c_2-t_1)\\
		0&b\xi^2&t_2&0&b\xi(t_2-c_1)&0\\
		\xi&0&0&\xi^2(t_1-c_2)/b&0&-t_1\\
		0&\xi&\xi^2(c_1-t_2)/b&0&-t_2&0
	\end{psmallmatrix}.
\end{align*}

Next, we check $[a_0]^2=\mu(\id-t^2)$. Let us set $A=[a_0]^2$. We get
\begin{enumerate}[leftmargin=*]
	\item $A_{11}=A_{22}=-2bf+c_1c_2=\mu(1-b^2)$.
	\item $A_{35}=A_{53}=A_{46}=A_{64}=0$ as they are scalar multiple of $bf-c_1t_1-c_2t_2+2t_1t_2=0$
	\item All other entries are zero.
\end{enumerate}

Now we deal with $B=[a_0][a_1]+[a_1][a_2]+[a_2][a_0]$. We get the following.
\begin{enumerate}[leftmargin=*]
	\item A clean computation shows that $B$ is a diagonal matrix.
	\item $B_{11}=B_{22}=3bf+3c_1c_2=\lambda(1-b^2)$. Indeed:
	\item In a similar fashion, we check that:
	\begin{align*}
		B_{33}=B_{44}&=3(-\xi^2 bf+\xi c_1c_2-\xi c_1t_1-\xi c_2t_2-\xi^2 t_1t_2)
		=\lambda(1-\xi b^2),\\
		B_{55}=B_{66}&=3(-\xi bf+\xi^2c_1c_2-\xi^2c_1t_1-\xi^2c_2t_2-\xi t_1t_2)
		=\lambda(1-\xi^2 b^2).
	\end{align*} 
\end{enumerate}
This shows that $N$ is an $\A$-module.

\medbreak

Now we show that $N$ is simple. Assume $N$ is not simple and let $0\neq N'\not\subseteq N$ be a submodule. In particular $N'_{|}=(S_j^+)^a\oplus(S_j^-)^b\oplus(M_j)^c$, where $0\leq a, b\leq 1, 0\leq c\leq 2$. Recall the basis $\lg x,y,v_1,v_2,w_1,w_2\rg$ of $N$.

If $a>0$ then $x\in N'$. We act with $a_0$ and get that $v_1\in N'$ (so $w_1\in N'$). Now, if $t_2\neq0$, then acting by $a_0$ again we have that $v_2\in N'$. Thus, $y\in N'$ and $N'=N$, a contradiction. On the other hand, if $t_2=0$, then we have $f_+=0$. Hence $g\neq 0$ and thus $c_1\neq 0$. So, $y\in N'$ and therefore $v_2\in N'$. Again $N'=N$. Thus, necessarily $a=0$.

Similarly, we show that $b=0$. Then $N_{|}\simeq M_j^c$, but there is no such $\A$-module when $3\mu\neq\lambda$.

For the converse, we start with the decomposition into $\G_{3,\ell}$-factors: $N$ cannot be a direct sum of modules of dimension 1, or modules of dimension 2, by the result in the previous section. And using Lemma \ref{lem:mixed-Mj-1} we discard:
\begin{align*}
	N_{|}&\simeq S_j^{+3}\oplus S_j^{-}\oplus M_j, & N_{|}&\simeq S_j^{+}\oplus S_j^{-3}\oplus M_j, & N_{|}&\simeq S_j^{+2}\oplus S_j^{-2}\oplus M_j.
\end{align*}
Hence \eqref{eqn:isotypic} holds. Notice that this is independent of the fact that $3\mu\neq \lambda$.

Now, as $3\mu\neq \lambda$, then $f\neq 0$. As $N$ is assumed to be simple, $a_0\cdot x\neq 0$. 
Now, $a_0\cdot x\in\lg y,v_1,v_2,w_1,w_2\rg$. If $a_0\cdot x\in\lg y\rg$, then $\lg x,y\rg$ defines a submodule. In other words, $\delta_{+}=(\delta_{+1} \, \delta_{+2})\neq 0$ in \eqref{eqn:a0matrix}. Let us assume that $\delta_{+1}\neq 0$ and rearrange the basis so that  $v'_1=\delta_{+1}v_1+\delta_{+2} v_2$, $w'_1=\delta_{+1}w_1+\delta_{+2} w_2$. In a similar way, we get that $\delta_{-}=(\delta_{-1} \, \delta_{-2})\neq 0$. If $\delta_{-2}=0$, then we get a contradiction to equations \eqref{eqn:general2}, as $f\neq 0$. This shows that we can assume $\delta_{+}=(1 \, 0)$, $\delta_{-}=(0 \, 1)$. The conditions on $\eta$ and $\sigma$ are now straightforward consequences from equations \eqref{eqn:general2}. The same applies for $\theta$ and $\eta$, which follow from equations \eqref{eqn:general3}, which become:
\begin{align*}
	&\begin{cases}
		b\tau_{11}+\theta_{11}=0\\
		b\tau_{11}-\theta_{11}=0
	\end{cases}
	&&&
	&\begin{cases}
		b\tau_{22}+\theta_{22}=0\\
		b\tau_{22}-\theta_{22}=0
	\end{cases}
	&&&
	&\begin{cases}
		b\tau_{12}+\theta_{12}=c_1\\
		b\tau_{21}-\theta_{21}=-c_2.
	\end{cases}
\end{align*}
Thus the result follows as this implies $\theta_{11}=\theta_{22}=0$,  $\tau_{11}=\tau_{22}=0$, and we set $t_1=\tau_{12}$, $t_2=\tau_{21}$. 
\epf

When $3\mu=\lambda$ we cannot define such a module. We observe that in this case $f=0$, $g=\mu(1-b^2)$ and $f_+=f_-=\mu$.
\begin{corollary}
	If $3\mu=\lambda$, then there are no 6-dimensional simple modules.
\end{corollary}
\pf
Assume $N$ is such a module. We follow the first lines of the proof of the converse in Proposition \ref{pro:dimension6} to see that \eqref{eqn:isotypic} holds. We write the action of $a_0$ as in \eqref{eqn:a0matrix} and conclude that we can assume $\delta_+=(1\,0)$. As well, as $g\neq 0$, we can set $\alpha=g$, $\beta=1$. It follows as well that $\delta_{-}\neq 0$. 

Now, as $3\mu=\lambda$, then $f=0$ and we cannot make any further assumptions on $\delta_{-}$. From $\delta_{+}\sigma=\delta_{+}\eta=0$ we get that $\sigma=\begin{psmallmatrix}
	0\\\sigma_2
\end{psmallmatrix}$, $\eta=\begin{psmallmatrix}
	0\\\eta_2
\end{psmallmatrix}$, with $\sigma_2\eta_2\neq 0$ as otherwise $\lg v_1,v_2,w_1,w_2\rg$ would define a submodule. Now we turn to $\delta_{-}\eta=\delta_{-}\sigma=0$ to see that $\delta_{-}=(d\,0)$ for some $d\neq 0$. 

We thus have: $a_0\cdot x=y-bv_1+w_1$, $a_0\cdot y=g\,x+dbv_1+dw_1$ and
\begin{align*}
	a_0\cdot v_1=\theta_{11}v_1+\theta_{12}v_2+\tau_{12}w_1+\tau_{12}w_2.
\end{align*}
Let us fix $\theta=\begin{psmallmatrix}
	\theta_{11}&\theta_{12}\\\theta_{21}&\theta_{22}
\end{psmallmatrix}$, $\tau=\begin{psmallmatrix}
	\tau_{11}&\tau_{12}\\\tau_{21}&\tau_{22}
\end{psmallmatrix}$, so \eqref{eqn:general3} becomes:
\begin{align*}
	\begin{cases}
		b\tau_{11} +\theta_{11}=d\\
		b\tau_{11} -\theta_{11}=-g/d,
	\end{cases}& 
	\begin{cases}
		b\tau_{12} +\theta_{12}=0,\\ 
		b\tau_{12} -\theta_{12}=0,\\ 
	\end{cases}&			
	\begin{cases}
		b\tau_{22}-\theta_{22}=g/d,\\
		b\tau_{22}+\theta_{22}=-d.
	\end{cases}			
\end{align*}
so $\tau_{12}=\theta_{12}=0$, which imply that $\lg v_1,w_1\rg$ is a submodule.
\epf

\subsubsection{Isomorphism classes of 6-dimensional simple modules}\label{sec:isoclasses}

In this part we compute the sets $\Up_j$, $j\in\J$, and choose in Definition \ref{def:6dim} representatives for isomorphism classes of the modules $N_j(\up)$, $\up\in\Up_j$.

We start with a simple observation.
\begin{remark}\label{rem:gfh} We fix $j\in\J$, $b=b_j$. 
	Set $f_{+}=f_{+j}$, $g=g_j$, $h=h_j$.
	\begin{itemize}[leftmargin=*]
		\item $g=0$ if and only if  $\mu=-2\frac{\lambda}{3}$; hence $h=\frac{\lambda}{3}(1+b^2)$ and $f_+=b^2\lambda/3$.
		\item If $f_+=0$, then $\mu=\frac{2+b^2}{b^2-1}\frac{\lambda}{3}$ and $h=\frac{\lambda}{3}$.
	\end{itemize}
	Thus we cannot have $f_+=0$ and $g=0$ simultaneously; same for $f_+$ and $h$. 
	On the other hand, $g=h=0$ only if $b^2=-1$ (hence $\ell$ is even).
\end{remark}

We follow with a technical lemma whose impact will become apparent in Proposition \ref{pro:dimension6} below. In preparation, we set $\sol_j$, $j\in\J$, as
\begin{align}\label{eqn:solution}
	\sol_j\coloneqq(h_j+\sqrt{h_j^2-4f_{+j}g_j})/{2}.
\end{align}

\begin{remark}\label{rem:c+1}
By Remark \ref{rem:gfh}, $\sol_j=0$ if and only if $g_j=h_j=0$, which implies $\mu=-2\frac{\lambda}{3}$ and $b_j^2=-1$. 
\end{remark}

\begin{lemma}\label{lem:solutions}
Let $\Up_j$ be as in Definition \ref{def:solutions}.
\begin{enumerate}[leftmargin=*]
\item[(i)] If $\mu=-2\frac{\lambda}{3}$ and $b_j^2=-1$, then $\Up_j\coloneqq\{(0,0,c,-\frac{\lambda}{3c}): c\in\k^\times\}$.
\item[(ii)] Otherwise, $\Up_j=\Up_j^{(1)}\cup \Up_j^{(2)}$, for
\begin{align*}
	\Up_j^{(1)}=\{(c,\frac{g_j}{c},\frac{\sol_j}{c},\frac{cf_{+j}}{\sol_j}) : c\in\k^\times\},  \	\Up_j^{(2)}=\{(\frac{g_j}{c},c,\frac{cf_{+j}}{\sol_j},\frac{\sol_j}{c}) : c\in\k^\times\}.
\end{align*}
\end{enumerate}
As well, 
\begin{align}\label{eqn:discriminant}
\Up_j^{(1)}=\Up_j^{(2)}\iff h_j^2=4f_{+j}g_j;
\end{align}
and $\Up_j^{(1)}\cap \Up_j^{(2)}=\emptyset$ otherwise.
\end{lemma}

We shall look deeper into condition \eqref{eqn:discriminant} after the proof, see \eqref{eqn:c} below.
\pf
We start with two facts:
	\begin{itemize}[leftmargin=*]
		\item $c_1=c_2=0\iff g=0$ ($\mu=-2\frac{\lambda}{3}$) and $h_j=0$ ($b_j^2=-1$, $\ell$ even). 
		\item We cannot have $t_1=t_2=0$, as this gives $f_+=h=0$. 
	\end{itemize}
Hence, the solutions to \eqref{eqn:conditions-six} in context (i), as $t_1t_2=f_+=-\frac{\lambda}{3}\neq 0$, are of the form $(0,0,t_1,f_+/t_1)\in\Up_1$.

Now we fix $j\in\J$ and study different alternatives to show $(ii)$. In what follows we omit subscripts referring to $j$.
	\begin{enumerate}[leftmargin=*]
		\item Assume $gf_+\neq 0$. In particular $c_1c_2=g\neq 0$ and $t_1t_2=f_+\neq 0$. So $t_1$ is part of a solution $\up=(c_1,c_2,t_1,t_2)$ if and only if $(c_1t_1)^2-h(c_1t_1)+gf_+=0$. Thus \[
		t_1=\sol^\pm(c_1)\coloneqq(h\pm\sqrt{h^2-4f_+g})/(2c_1), \qquad t_2=f_+/t_1.
		\]
We see:
		\[
		\up=(c_1,g/c_1,\sol^+(c_1),f_+/\sol^+(c_1))=(c,g/c,\sol/c,cf_+/\sol)\in\Up^{(1)},
		\]
		for $c=c_1$ or, as $\sol^-(c_1)=f_+g/\sol^+(c_1)$ and thus, for $c=g/c_1$:
		\[
		\up=(c_1,g/c_1,\sol^-(c_1),f_+/\sol^-(c_1))=(g/c,c,cf_+/\sol,\sol/c)\in\Up^{(2)}. 
		\]
		\item If $g=0$, $j>1$, then the solutions are either $\up=(c_1,0,h/c_1,c_1f_+/h)$ or $\up=(0,c_1,c_2f_+/h,h/c_2)$, both of which belong to $\Up^{(1)}\cup\Up^{(2)}$ as $\sol=h$.
		\item If $f_+=0$, then $c_1c_2=g\neq 0$ and thus the solutions are either $\up=(c_1,g/c_1,h/c_1,c_1f_+/h)$ or $\up=(g/c_2,c_2,c_2f_+/h,h/c_2)$, again in $\Up^{(1)}\cup\Up^{(2)}$ as $\sol=h$.
\end{enumerate}
Finally, let us set $\sol^\pm=\sol^\pm(1)$, so  $\sol=\sol^+$. Thus $\up=(c,g/c,\sol/c,cf_+/\sol)=(g/c',c',c'f_+/\sol,\sol/c')\in \Up^{(1)}\cap\Up^{(2)}$ if and only if $\sol^2=gf_+=\sol^+\sol^-$, if and only if $\sol^+=\sol^-$: that is $h^2-4f_+g=0$. Conversely, if $h^2-4f_+g$ then $\Up^{(1)}=\Up^{(2)}$.
\epf

We look into condition \eqref{eqn:discriminant} from Lemma \ref{lem:solutions}. On the one hand:
\begin{align}\label{eqn:h2}
	h_j^2-4f_{+j}g_j=-\frac{(b_j^2-1)^2}{3}\left(\mu^2+\frac{2}{3}\lambda\mu-\frac{\lambda^2}{3}\frac{(b_j^2+1)^2}{(b_j^2-1)^2}\right).
\end{align}
Hence, $h_j^2=4f_{+j}g_j$ if and only if $\mu=\c_{\pm j}\frac{\lambda}{3}$, for
	\begin{align}\label{eqn:c}
		\c_{\pm j}=\frac{1-b_j^2 \pm 2\sqrt{b_j^4+b_j^2+1}}{b_j^2-1}.
	\end{align}
Thus $\c_{+ j}\c_{- j}=-3(\frac{b_j^2+1}{b_j^2-1})^2$. We collect some other useful observations next.
\begin{remark}\label{rem:c}
Keep in mind that $\c_{\pm j}$ depends on $\ell$. 
\begin{itemize}[leftmargin=*]
\item By \eqref{eqn:h2}, $\c_{\pm j}=\c_k$ if and only if $	j=k$ in $\J$. 
\item If $b_j^2=-1$, then $\c_{+j}=-2$; we discard $\c_{-j}=0$ as $(\lambda,\mu)\neq (0,0)$.
\item It is easy to check that $\c_{\pm j}\neq 1$ for any $\ell$ and any $j\in\J$. 
\item Notice that $\c_{\pm1}=\c_{\pm2}=-1$ when $\ell=3$.
\end{itemize}
\end{remark}

Following the description of $\Up_j$ in Lemma \ref{lem:solutions}, we define particular modules to characterize isomorphism classes in Corollary \ref{cor:isoclasses}.
\begin{definition}\label{def:6dim}
Assume $3\mu\neq\lambda$. 
\begin{itemize}[leftmargin=*]
\item If $g=0$ and $b_j^2=-1$ (so $\mu=\c_{+j}\frac{\lambda}{3}=-2\frac{\lambda}{3}$), then  we set
\begin{align}\label{eqn:Nj-g0j1}
	N_j^{(1)}&\coloneqq N_j(0,0,1,-\lambda/3).
\end{align}
\item Otherwise, we define, for each $j\in\J$:
\begin{align}\label{eqn:Nj-def}
	N_j^{(1)}&\coloneqq N_j(1,g,\sol,f_+/\sol), & 	N_j^{(2)}&\coloneqq N_j(g,1,f_+/\sol,\sol).
\end{align} 
\end{itemize}
\end{definition}

\begin{corollary}\label{cor:isoclasses}
Let $N=N_j(\up)$, $\up\in\Up_j$.
\begin{enumerate}
\item If $\mu=-2\frac{\lambda}{3}$ and $b_j^2=-1$ (so $\c_{+ j}=-2$), then $N\simeq N_j^{(1)}$ as in \eqref{eqn:Nj-g0j1}.
\item Otherwise, $N\simeq N_j^{(1)}$ or $N\simeq N_j^{(2)}$ as in \eqref{eqn:Nj-def}; and $N_j^{(1)}\simeq N_j^{(2)}$ if and only if $h_j^2-4f_{+j}g_j= 0$ (so $\mu=\c_{\pm j}\frac{\lambda}{3}$).
\end{enumerate}
In other words, there is a unique, up to isomorphism, 6-dimensional simple module $N$ if and only if $\mu=\c_{\pm j}\frac{\lambda}{3}$, in which case $N\simeq N_j^{(1)}$; this happens at most for a single $j\in\J$.
\end{corollary}
\pf
Follows by Lemma \ref{lem:solutions}, by normalizing $c_1=1$ or $c_2=1$. The last part follows from Lemma \ref{lem:solutions}, combined with Remark \ref{rem:c}.
\epf
%

\subsection{Classification of simple modules}\label{sec:classification}

In this section we present the collection of all simple $\mA_{\lambda,\mu}$-simple modules, up to isomorphism. 

Let us recall Remark \ref{rem:dim1-notprojective} which explains that the modules $L_0^\pm$ are not projective. As well, in Section \ref{sec:extensions} we shall compute the spaces $\Ext^1(L,L')$ for every pair of simple modules $L,L'$. In particular, when these computations are focused on the simple $\mA_{\lambda,\mu}$-modules of dimension two $T_j$ and  $L_j(a,c)$ defined so far, they readily provide the next result, which will become a useful technical tool to round up the classification. More explicitly, we compute non-trivial extensions in $\Ext^1(L_0^\pm,L_0(a,\pm a))$ in Lemma \ref{lem:extSM} and in $\Ext^1(L_j(a,c), T_j)$, $j\in\J$, in Lemma \ref{lem:extTM} (when $\mu=\frac\lambda3$).

We extract here an immediate consequence of these computations.  

\begin{lemma}\label{lem:1and2-not-projective}
If $L$ is a simple $\mA_{\lambda,\mu}$-module of dimension 2, then $L$ is not projective. In particular $\dim P(L)\geq 3$.\qed
\end{lemma}

As well, we need to anticipate some ideas on projective covers of simple modules, which we will work on full scope on \S \ref{sec:projective}. 

\begin{proposition}\label{pro:projective0}
$\mA(S_0^\pm)$ is the projective cover of the $\mA$-module $L_0^\pm$.
\end{proposition}
\pf
We proceed as in \cite[Proposition 5.4 (ii)]{GI}. We analyze the case $S_0^+$, the case $S_0^-$ being analogous.
It suffices to show that it is indecomposable. Assume $\mA(S_0^+)=U_1\oplus U_2$, we can assume that $S_0^+\subset (U_1)_{|}$. Let $x$ be such that $S_0^+=\lg x\rg$: there are $c_1, c_2\in \k$ with $x=c_1x_1+c_2x_2\in U_1$, where $x_1=1\ot z$, $x_2=(a_0a_2a_0a_1-\lambda  a_0a_1)\ot z$. Furthermore, we can assume $c_1\neq 0, c_2\neq 0$.

Notice that $a_0\cdot x_2=0$, so $a_0\cdot x=c_1a_0\ot z$ y $a_0\ot z\in U_1$.
If we act by $s$, $a_0$, $a_2$ and by $a_0$, in that order, we get:
\[
a_0\ot z \stackrel{s}{\to}  a_1 \ot z \stackrel{a_0}{\to}  a_0a_1\ot z \stackrel{a_2}{\to}  a_2a_0a_1\ot z \stackrel{a_0}{\to}  a_0a_2a_0a_1\ot z\in U_1.
\]
As $U_1$ is a submodule, $x_2\in U_1$. Thus  $x_1=1\ot z\in U_1$ and hence $I(S_0^+)\subset U_1$, $U_2=\{0\}$. Therefore this module is indecomposable. 
\epf

\begin{remark}
It can be shown that 	
\[
\mA(S_j^\pm)\simeq S_j^{+2}\oplus S_j^{-2}\oplus M_j^{4}.
\]
We choose to spare the reader this (lengthy) computation here, see \cite{tesis}.
\end{remark}

These results allows us to state the following partial result.

\begin{lemma}\label{lem:dim7}
There are no simple $\mA_{\lambda,\mu}$-modules of dimension $\geq 7$.
\end{lemma}
\pf
Assume $L$ is such a module, and let $j\in\J_0$ be such that $L=L[j]$. Let $P(L)$ be its projective cover. In particular,
\[
\dim L\dim P(L)\geq 49.
\]

If $j=0$ and $\lambda\neq 0$, as $\A(S_0^\pm)$ is the projective cover of the modules $L_0^\pm$, we obtain the computation
\begin{align*}
72&\geq 1\cdot 12+1\cdot 12+4\cdot 2\cdot \dim P(M_0(a,c))+\dim L\dim P(L)\\
&\geq 40+\dim L\dim P(L).
\end{align*}
which is a contradiction. If $\lambda=0$, the computation is
\begin{align*}
	72&\geq 1\cdot 12+1\cdot 12+ 2\cdot \dim P(M_0)+\dim L\dim P(L)\\
	&\geq 30+\dim L\dim P(L),
\end{align*}
which yields a contradiction as well.

Assume now that $j\in\J$. If $3\mu=\lambda$ we get
\begin{align*}
	72&\geq 2\cdot \dim P(T_j)+4\cdot 2\cdot \dim P(L_j(a,c))+\dim L\dim P(L)\\
	&\geq 30+\dim L\dim P(L).
\end{align*}
which is again a contradiction.
On the other hand, if $3\mu=\lambda$, we can have either two modules of dimension 6, say $N_j^{(1)}$ and $N_j^{(2)}$ or a unique one, $N_j$. In the first case, we get
\begin{align}
\begin{split}\label{eqn:formulaNj-1}
72&\geq 6\cdot P(N_j^{(1)})+6\cdot P(N_j^{(1)})+\dim L\dim P(L)\\
&\geq 72+\dim L\dim P(L).
\end{split}
\end{align}
In the second, we get
\begin{align}\label{eqn:formulaNj-2}
72\geq 6\cdot P(N_j)+\dim L\dim P(L)\geq 36+\dim L\dim P(L).
\end{align}
In any case, we get a contradiction. 
\epf

Recall the definition of $\c_{\pm j}\in\k$ in \eqref{eqn:c}.
\begin{corollary}\label{cor:projective-Nj}
If 	$\mu=\c_{\pm j}\frac{\lambda}{3}$ then the projective cover of the unique 6-dimensional $\mA_{\lambda,\mu}$-module $N_j^{(1)}$ is the induced module $\mA(S_j^+)$. 
Otherwise, the two simple $\mA_{\lambda,\mu}$-modules $N_j^{(1)}$ and $N_j^{(2)}$ of dimension six are projective.
\end{corollary}
\pf
Follows from the proof of Lemma \ref{lem:dim7}, more precisely the first part follows from \eqref{eqn:formulaNj-1}, which gives $\dim P(N_j)=6=\dim N_j$, and thus $N_j=P(N_j)$. The second assertion follows from \eqref{eqn:formulaNj-2}, which implies $\dim P(N_j)=12$. As $\mA(S_j^+)$ is projective and $N_j$ is the unique simple module in this setting, we get that $\mA(S_j^+)$ projects onto $N_j$, and thus onto the projective cover $P(N_j)$. This shows $P(N_j)=\mA(S_j^+)$ by a dimension argument.
\epf


\begin{theorem}\label{thm:simple}
	The simple modules for $\mA_{\lambda,\mu}$ are, up to isomorphism:
	\begin{enumerate}[leftmargin=*]
		\item If $(\lambda,\mu)=(0,0)$, then there are $3\ell$ modules $L_S$, with $L_{S|}\simeq S$, $a_0=0$, one for each $S\in\widehat{\G_{3,\ell}}=\{S_j^+,S_j^-,M_j:j\in\J_0\}$.
		\item 
		If $\mu=\frac\lambda3\neq 0$, then there are $5\ell+1$  modules, namely:
		\begin{itemize}[leftmargin=*]
			\item Two 1-dimensional modules $L_0^+$ and $L_0^-$, so that $L_{0|}^\pm=S_0^\pm$ and $a_0=0$.
			\item Four 2-dimensional modules  $L_j(a,c)$, $a,c\in\{\pm\sqrt{\lambda/3}\}$ for each $j\in\J_0$:
			\begin{align}\label{eqn:teo-2dim}
			L_j(a,c)_{|}&\simeq M_j, &	[a_0]&=\begin{pmatrix} a&-\zeta^jc\\c&-a \end{pmatrix}, a,c\in\{\pm\sqrt{\lambda/3}\}.
			\end{align}
			\item One 2-dimensional module $T_{j}$, for each $j\in\J$, with
			\begin{align*}
			T_{j|}&\simeq S_j^+\oplus S_j^-, &	[a_0]&=\begin{pmatrix} &\mu(1-\zeta^j)\\1&0 \end{pmatrix}.
		\end{align*}
		\end{itemize}
		\item 
		If $\mu\neq \dfrac\lambda3$ and $\lambda\neq 0$, then we distinguish three cases:
		\begin{enumerate}[leftmargin=*]
			\item If $\lambda\neq0$, $\mu\neq\c_{\pm j}\frac{\lambda}{3}$, any $j\in\J$, then there are $2\ell+4$  modules:
			\begin{itemize}[leftmargin=*]
				\item Two 1-dimensional modules $L_0^\pm$, with $a_0=0$.
				\item Four 2-dimensional modules $L_0(a,c)$  as in \eqref{eqn:teo-2dim}--so $[a_0]=\begin{psmallmatrix}
					a&-c\\c&-a
				\end{psmallmatrix}$.
				\item Two 6-dimensional modules $N_j^{(1)}$ and $N_j^{(2)}$ as in \eqref{eqn:Nj-def}, for $j\in\J$.
			\end{itemize}
			\item If $\mu=\c_{\pm j}\frac{\lambda}{3}\neq 0$, then there are $2\ell+3$  modules, namely:
			\begin{itemize}[leftmargin=*]
				\item Two 1-dimensional modules $L_0^\pm$.
				\item Four 2-dimensional modules $L_0(a,c)$.
				\item One 6-dimensional module $N_j^{(1)}$ as in \eqref{eqn:Nj-g0j1} and two 6-dimensional modules $N_k^{(1)}$ and $N_k^{(2)}$ for each $k\in\J$, $k\neq j$, as in \eqref{eqn:Nj-def}.
			\end{itemize}
		\item If $\lambda=0$, $\mu\neq0$, then there are $2\ell+1$  modules, namely:
		\begin{itemize}[leftmargin=*]
			\item Two 1-dimensional modules $L_0^\pm$.
			\item One 2-dimensional  module $L_0(0,0)$ as in \eqref{eqn:teo-2dim}--so $[a_0]=\begin{psmallmatrix}
				0&0\\0&0
			\end{psmallmatrix}$.
			\item Two 6-dimensional modules $N_j^{(1)}$ and $N_j^{(2)}$ as in \eqref{eqn:Nj-def}, for $j\in\J$.
					\end{itemize}
		\end{enumerate}
	\end{enumerate}
\end{theorem}
\pf
The graded case is Proposition \ref{pro:graded-case}. For the rest, we have listed all possible simple modules of dimension 1, 2 and 6; and shown that there are no such modules of dimensions 3, 4, 5 or $\geq 7$. The result follows.
\epf

\subsection{On the tensor structure of $\Rep \mA_{\lambda,\mu}$}\label{sec:tensor-structure}

\

Let $L=L_j^{p,q,r}(\alpha,\beta,\delta_{\pm},\eta,\sigma,\theta,\tau)$ be an indecomposable $\mA_{\lambda,\mu}$-module, as in Lemma \ref{lem:general-a0}. That is, $L_{|}\simeq S_j^{+p}\oplus S_j^{-q}\oplus M_j^{r}$ and the module structure is codified by a list of matrices $\alpha,\beta,\delta_{\pm},\eta,\sigma,\theta,\tau$ as in \eqref{eqn:a0matrix}.

In particular, any $\mA_{\lambda,\mu}$-module $M\in\Rep \mA_{\lambda,\mu}$ is a direct sum of a collection of such modules. We start by computing their isomorphism classes in Lemma \ref{lem:repA-isos}.
Then, we describe the tensor structure of this category. Namely we compute in Lemma \ref{lem:repA-duals} the duals $L_j^{p,q,r}(\alpha,\beta,\delta_{\pm},\eta,\sigma,\theta,\tau)^*$ and we write down in  Proposition \ref{pro:tensor} the tensor products 
\[L_j^{p,q,r}(\alpha,\beta,\delta_{\pm},\eta,\sigma,\theta,\tau)\ot L_k^{p',q',r'}(\alpha',\beta',\delta_{\pm}',\eta',\sigma',\theta',\tau').\]

We shall write $L_{-j}$ and $L_{j+k}$ to indicate the subscripts corresponding to the values $i\in\{0,\dots,\ell-1\}$ so that $i\equiv -j\mod\ell$ and $i\equiv j+k\mod\ell$.

\begin{lemma}\label{lem:repA-isos} There is an isomorphism of $\mA_{\lambda,\mu}$-modules
\[
L_j^{p,q,r}(\alpha,\beta,\delta_{\pm},\eta,\sigma,\theta,\tau)\simeq  L_k^{p',q',r'}(\alpha',\beta',\delta_{\pm}',\eta',\sigma',\theta',\tau')
\] if and only if $(j,p,q,r)=(k,p',q',r')$ and 
\begin{align*}
\alpha'&=X\beta Y^{-1}, & \beta'&=Y\beta X^{-1}, & \delta_+'&=V\delta_+ X^{-1}, & \delta_-'&=V\delta_- Y^{-1},\\
\eta'&=X\eta V^{-1}, & \sigma'&=Y\sigma V^{-1}, &\theta'&=V\theta V^{-1}, & \tau'&=V\tau V^{-1}.
\end{align*} 
for some invertible matrices $X\in\k^{p\times p}$, $Y\in\k^{q\times q}$, $V\in\k^{r\times r}$.
\qed
\end{lemma}

\begin{lemma}\label{lem:repA-duals}
If $L=L_j^{p,q,r}(\alpha,\beta,\delta_{\pm},\eta,\sigma,\theta,\tau)$, then
\[
L^\ast \simeq L_{-j}^{p,q,r}(-b^{-1}\beta,-b^{-1}\alpha,-\eta,-\sigma,\delta_{\pm},\tau,-\theta).
\]
\end{lemma}
\pf
We have already established that $(L^\ast)_{|}\simeq S_{\ell-j}^{+p}\oplus S_{\ell-j}^{-q}\oplus M_{\ell-j}^r$ and, from $S(a_0)=-t^{-1}a_0$ we get that the matrix $[a_0^\ast]$ of the action of $a_0$ in the dual basis $\{\delta_x,\delta_y,\delta_v,\delta_w\}$ of the basis $\{x,y,v,w\}$ of $L$ is given by $[a_0^\ast]=-[t^{-1}]^t[a_0]^t$, that is:
\begin{align*}
[a_0^\ast]=&-\begin{psmallmatrix}
	b^{-1}&0&0&0\\
	0&-b^{-1}&0&0\\
	0&0&0&b^{-2}\\
	0&0&1&0\\
\end{psmallmatrix}\begin{psmallmatrix}
0&  \beta&  -b  \delta_+&\delta_+\\
\alpha&0& b  \delta_- &\delta_-\\
\eta&\sigma&  \theta&\tau\\
-b   \eta&   b  \sigma&-\zeta^j  \tau &-  \theta
\end{psmallmatrix}=\begin{psmallmatrix}
	0&-b^{-1}\beta&\delta_+&-b^{-1}\delta_+\\
	-b^{-1}\alpha&0&\delta_-&b^{-1}\delta_-\\
	b^{-1}\eta & -b^{-1}\sigma& \tau& b^{-2}\theta\\
	-\eta & -\sigma &  -\theta   &-\tau 
\end{psmallmatrix}.
\end{align*}
The lemma follows.
\epf

\begin{proposition}\label{pro:tensor}
Let $L_j,L_k\in\Rep\mA_{\lambda,\mu}$ denote the modules
\[
L_j=L_j^{p,q,r}(\alpha,\beta,\delta_{\pm},\eta,\sigma,\theta,\tau), \ L_k=L_k^{p',q',r'}(\alpha',\beta',\delta_{\pm}',\eta',\sigma',\theta',\tau').
\] 
Then the following isomorphism holds:
	\begin{align*}
L_j\ot L_k
\simeq L_{j+k}^{{\bf p},{\bf q},{\bf r}}(\alpha^{\ot},\beta^{\ot},\delta^{\ot}_{\pm},\eta^{\ot},\sigma^{\ot},\theta^{\ot},\tau^{\ot}).
	\end{align*}
%
where the integers ${\bf p}, {\bf q}, {\bf r}$ are as follows
\begin{align*}
	{\bf p}&=pp'+qq'+rr', &
	{\bf q}&=pq'+qp'+rr', &
	{\bf r}&=2(pr'+rp'+rq'+rr').
\end{align*}
Set $b=b_jb_k$, then the parameter matrices are given by
\begin{align*}
	\alpha^\ot&= 
	\begin{bsmallmatrix} 
		b_j\balpha'&\bbeta&0\\ 
		\balpha&-b_j\bbeta'&0\\
		0&0&\btheta-b\btau'
	\end{bsmallmatrix}, &&
	\delta_+^\ot=
\begin{bsmallmatrix} 
	\bdelta_+'&0& b_k^{-1}\bdelta_+&0&0 \\ 
	0&\bdelta_-'&0&-b_k^{-1}\bdelta_-&0\\ 
	\beeta & -\bsigma & -b_j\beeta'& -b_j\bsigma'& b_k^{-1}(\btau-\btheta') 
\end{bsmallmatrix},
\end{align*}
\begin{align*}
\beta^\ot=
\begin{bsmallmatrix} 
	b_j\bbeta'&\bbeta&0\\
	\balpha &-b_j\balpha'&0\\
	0&0&\btheta+b\btau'
\end{bsmallmatrix},
&&
	\delta_-^\ot=
	\begin{bsmallmatrix} 
		\bdelta'_- & 0 & 0 & - b_k^{-1}\bdelta_+ & 0\\
		0& \bdelta'_+ & b_k^{-1}\bdelta_-  & 0 & 0 \\
		\beeta & -\bsigma & -b_j \beeta' & - b_j \bsigma' & b_k^{-1}(\btau -\btheta')
	\end{bsmallmatrix},
\end{align*}
\begin{align*}
	\eta^\ot& = 
	\begin{bsmallmatrix} 
		b_j\beeta' & 0  & \frac{b_k^{-1}}{2}\bdelta_+  \\
		0& -b_j\bsigma' & \frac{b_k^{-1}}{2}\bdelta_-  \\
		\beeta  & 0 & - \frac12 \bdelta'_+   \\
		0& \bsigma & \frac12 \bdelta'_-  \\
		0& 0 &-\frac{b_j}{2}(\btau+\btheta')\
	\end{bsmallmatrix},
	&
	\sigma^\ot&=
	\begin{bsmallmatrix} 
		b_j\bsigma'   & 0 & -\frac{b_k^{-1}}{2}\bdelta_+\\
		0& - b_j\beeta' &- \frac{b_k^{-1}}{2}\bdelta_-  \\
		0 &  \bsigma  & \frac12 \bdelta'_+\\
		\beeta & 0 &  -\frac12 \bdelta'_-  \\
		0 & 0 & -\frac{b_j}{2}(\btau+\btheta')
	\end{bsmallmatrix},
\end{align*}
\begin{align*}
	\theta^\ot&= 
	\begin{bsmallmatrix} 
		b_j\btheta'  &  \bbeta  & 0 & 0& 0\\
		\balpha & - b_j\btheta' &  0 & 0& 0\\
		0&0&\btheta  & 0 & \bdelta'_+  \\
		0&0&0& \btheta &  \bdelta'_-  \\
		0&0&-b_jb\beeta' & b_jb\bsigma' &  -\btheta   
	\end{bsmallmatrix}, 
&	\tau^\ot&=
	\begin{bsmallmatrix} 
		\btau' & 0 & 0 & 0 &-   b_k^{-1}b \bdelta_+\\
		0& \btau' & 0 & 0 & b_k^{-1}b\bdelta_- \\
		0 & 0 & b_k^{-1}\btau & -b_k^{-1}\bbeta' & 0\\
		0 & 0 & b_k^{-1}\balpha' & -b_k^{-1}\btau & 0\\
		-\bbeta &- \bsigma & 0 &0 &- \btau' 
	\end{bsmallmatrix}.
\end{align*}
For each $m\in\{\alpha, \beta, \delta_\pm, \eta, \sigma, \tau, \theta\}$ and $m'\in\{\alpha', \beta', \delta_\pm', \eta', \sigma', \tau', \theta'\}$ we denote
\begin{align*}
	\bm&=m\ot \id, & \bm'&=\id\ot\, m'
\end{align*}
for the appropriate size of the matrix $\id$ at each side of $\ot$. 
\end{proposition}
\pf
By Proposition \ref{pro:tensor-group}, if $\lg x, y, v, w\rg$ and $\lg x',y',v',w'\rg$ are basis of $L_j$ and $L_k$, respectively, then
$(L_j\ot L_k)_{|}\simeq S_{j+k}^{+{\bf p}}\oplus S_{j+k}^{-{\bf q}}\oplus M_{j+k}^{{\bf r}}$, for ${\bf p}, {\bf q}, {\bf r}$ as stated, 
in the basis 
\begin{align*}
	\{\bx_1&=x\ot x', \bx_2=y\ot y', \bx_3=b_j v\ot w'+b_k w\ot v', \\
	&	\by_1=x\ot y', \by_2=y\ot x', \by_3=b_j v\ot w'-b_k w\ot v',\\
	&\bv_1=x\ot v', \bv_2=y\ot v', \bv_3=v\ot x', \bv_4=v\ot y', \bv_5=w\ot w',\\
	& \bw_1=b_j x\ot w', \bw_2=-b_j y\ot w', \bw_3=b_kw\ot x', \bw_4=-b_kw\ot y', \\
	&\hspace*{8cm}\bw_5=(b_jb_k)^2 v\ot v'\}.
\end{align*}
Finally we can compute the action of $a_0$ on each element, for instance:
\begin{align*}
	a_0\cdot \bx_1&=a_0\cdot x\ot 1\cdot x' + t \cdot x \ot a_0\cdot x'\\
	&= (\beta y - b_j \delta_+v+\delta_+w)\ot x'+b_j x\ot (\beta' y' - b_k \delta'_+v'+\delta'_+w')\\
	&=b_j \bbeta'  \by_1 + \bbeta \by_2- b_jb_k \bdelta'_+ \bv_1 - b_j \bdelta_+\bv_3  + \bdelta'_+ \bw_1 +b_k^{-1}\bdelta_+ \bw_3.
\end{align*}
The result follows.
\epf

When $\ell=1$, so $\G_{3,\ell}\simeq \s_3$, we observe that we can assume $\mu=0$ (as $t^2=1$). We thus obtain a family of pointed Hopf algebras over $\s_3$, namely $\mH_\lambda=\mA_{\lambda,0}$, $\lambda\in\k$. We can further normalize $\lambda=1$ for the non-graded case. 

This is indeed how the unique, up to isomorphism, finite-dimensional non semisimple and non-graded pointed Hopf algebra $\mH_1\coloneqq \mA_{0,1}$ arises. When $\lambda=0$ we obtain the graded Hopf algebra $\mH_0=\fk_3\#\k\s_3$. The corresponding category $\Rep\mH_\lambda$ was studied in \cite{GI}. The results in this article apply, or rather restrict to (as well as are inspired by) this particular case. The computations in this subsection lead us to the following corollary. 

Recall that that $C_\ell$ stands for the cyclic group of order $\ell$, that in this part we shall represent by the quotient $\Z/\ell\Z$, which in turn we identify as a set with $\{0,\dots,\ell-1\}\subset \Z_{\geq0}$.

\begin{corollary}\label{cor:tensor}
Let $\R_{\lambda,\mu}=\Rep \mA_{\lambda,\mu}$ be the category of finite-dimensional representations of the Hopf algebra $\mA_{\lambda,\mu}$. 
Then $\R_{\lambda,\mu}=\bigoplus\limits_{j\in C_\ell}\R[j]$ is a $C_\ell$-graded tensor category, where 
$L$ is an object in $\R[j]$ if and only if $L=L[j]$. This grading is faithful and non-trivial if $\ell>1$.

Moreover, the tensor subcategory $\R[0]$ given by the trivial component is tensor equivalent to $\C_{\lambda}\coloneqq \Rep \mH_\lambda$. That is, $\R_{\lambda,\mu}$ is an extension of $\C_{\lambda}$.
\end{corollary}
\pf
The first part follows from Proposition \ref{pro:j-unico}, together with the results in this section, namely Lemma \ref{lem:repA-duals} and Proposition \ref{pro:tensor}. 

Next we show that $\R[0]\simeq_{\ot} \Rep \mH_\lambda$. To do this, let us fix a presentation of $\mH_\lambda$: this the algebra generated by elements $b_{(12)}$,  $b_{(13)}$, $b_{(23)}$ and $H_\tau$, $\tau\in\s_3$ so that $H_\tau b_{\sigma}=-b_{\tau\sigma \tau^{-1}}H_\tau$ and
\begin{align*}
	b_{(12)}^2&=0, & b_{(12)}b_{(23)}+b_{(23)}b_{(13)}+b_{(13)}b_{(12)}&=\lambda(1-H_{(123)}).
\end{align*}

The claim follows from the fact that for $V\in\R[0]$, the action $[t]$ of $t\in\G_{3,\ell}$ is such that $[t]^2=\id$, cf.~Lemma \ref{lem:uniquej}. 
Hence the defining relations \eqref{eqn:rels_def} for $\mA_{\lambda,\mu}$ are represented as
\begin{align*}
	[a_{(0)}]^2&=0, & [a_{(0)}a_{(1)}+a_{(1)}a_{(2)}+a_{(2)}a_{(0)}]&=\lambda(\id-[s]).
\end{align*}
We can thus define a functor $F\colon \R[0]\longrightarrow \Rep \mH_\lambda$, that is the identity both on morphisms as on the underlying vector spaces of $V\in\R[0]$, which in turn become modules over $\mH_\lambda$ by setting
\begin{align*}
[H_{(12)}]&=[t], & [H_{(123)}]&=[s], & [b_{(12)}]&=[a_0].
\end{align*}
By reversing this identification, we define a functor $G\colon \Rep \mH_\lambda\longrightarrow \R[0]$, which is an inverse for $F$, completing the proof.
\epf

\section{Projective covers}\label{sec:projective}

In this part we find the projective cover of each simple $\mA_{\lambda,\mu}$-module. Recall that we have already established in Proposition \ref{pro:projective0}
that $\mA(S_0^\pm)$ is the projective cover of the $\mA_{\lambda,\mu}$-module $L_0^\pm$. Recall the notation from \S\ref{sec:notation}.

We start by an observation on the modules $\mA(S_j^\pm)$, $j\in\J$.

\begin{lemma}\label{lem:iso-induced}
	Assume $g\neq 0$. Then $\mA(S_j^+)\simeq \mA(S_j^-)$.
\end{lemma}
\pf
Let us fix  $S_j^+=\lg z\rg\subset \mA(S_j^+)=A\ot \lg z\rg$, $S_j^-=\lg w\rg$.
Consider the isomorphism of $G$-modules $\varphi_0\colon \lg z\rg\to \lg (a_0+a_1+a_2)w\rg\subset \mA(S_j^-)$, given by $z\mapsto (a_0+a_1+a_2)w$. 
This naturally extends to $\varphi\colon \mA(S_j^+)\to \mA(S_j^-)$, by setting $a\, z\mapsto a(a_0+a_1+a_2)w$, $a\in \mA$.
We observe that
\[
\varphi((a_0+a_1+a_2)z)=(a_0+a_1+a_2)^2w=(3\mu+2\lambda)(1-\zeta^j)w=g\,w
\]
and thus $\varphi$ is an epimorphism -hence an isomorphism- when $g\neq 0$. 
\epf

\subsection{The projective cover of $L_j(a,c)$}\label{sec:claim}

In order to deal with the projective cover of the modules $L_j(a,c)$ (where $0<j\in\J$ only if $3\mu=\lambda$), we need to study the extensions of $L_j(a,c)$ by the tensor product $L_j(a,c)\ot L_j(a,c)$. This is inspired in \cite[Proposition 5.15]{GI}. Our goal is to prove the following.

\begin{claim}
For each $j\in\J_0$, $\dim\ext^1(L_j(a,c),L_j(a,c)\ot L_j(a,c))=1$. If we pick an indecomposable generator $P_j(a,c)\in \ext^1(L_j(a,c),L_j(a,c)\ot L_j(a,c))$, then this a projective cover of $L_j(a,c)$.\qed
\end{claim}

Recall that $L_j(a,c)=\lg v,w\rg$,  where
$[a_0]=\begin{psmallmatrix*} a&-\zeta^jc\\c&-a \end{psmallmatrix*}$,  $a,c\in\{\pm\sqrt{\lambda/3}\}$.

\begin{remark}
	Generically, there are {\it four} simple modules, for each $j\in\J_0$, non isomorphic of dimension 2. Some remarks:
	\begin{enumerate}
		\item If $3\mu\neq \lambda$, only exist for $j=0$.
		\item If $\lambda=0$, then $a=c=0$ (there exists an unique module, with $j=0$).
	\end{enumerate}
\end{remark}

We take a copy of $L_j(a,c)=\lg v_1,w_1\rg$ and one for $L_0(a,c)=\lg v_2,w_2\rg$, and
we shall to study the module $L^2=L_j(a,c)\ot L_0(a,c)$, with basis (canonical):
\[
\{v_1\ot v_2, v_1\ot w_2, w_1\ot v_2, w_1\ot w_2\}.
\]
It follows that, in this basis, $[s]=\begin{psmallmatrix}
		\xi^2&0&0&0\\0&1&0&0\\0&0&1&0\\ 0&0&0&\xi
	\end{psmallmatrix}$, $[t]=\begin{psmallmatrix}
		0&0&0&b^2\\0&0&b^2&0\\0&1&0&0\\ 1&0&0&0
	\end{psmallmatrix}$.
As well, as $\Delta(a_0)=a_0\ot 1+t\ot a_0$ we obtain that
\begin{align*}
	[a_0]=\begin{psmallmatrix*}
		a&0&b^2(a-c)&-b^2c\\0&a&b^2c&-b^2(a+c)\\a+c&-c&-a&0\\ c&c-a&0&-a
	\end{psmallmatrix*}.
\end{align*}
We next shift to a more suitable basis.
\begin{lemma}
There is a basis $\{x,y,v,w\}$ of $L=L_j(a,c)\ot L_0(a,c)$ so that $\lg x\rg_{|}\simeq S_j^+$, $\lg y\rg_{|}\simeq S_j^-$, $\lg v,w\rg_{|}\simeq M_j$  and thus
\[
L_{|}\simeq S_j^+\oplus S_j^-\oplus M_j.
\]
In particular, $[s]=\diag(1,1,\xi,\xi^2)$ and $[t]=\begin{psmallmatrix}
	b&0&0&0\\0&-b&0&0\\0&0&0&b^2\\ 0&0&1&0
\end{psmallmatrix}$. In this basis, the action of $a_0$ is determined by	
\begin{align*}
	[a_0]=\begin{psmallmatrix*}
		0&(bc-a)&-\frac{b(a+c)}{2}&\frac{b^2(a+c)}{2}\\
		-(a+bc)&0&\frac{b(a+c)}{2}&\frac{b^2(a+c)}{2}\\
		-b(a-c)&b(a-c)&-a&b^2c\\
		(a-c)&(a-c)&-c&a
	\end{psmallmatrix*}.
\end{align*}
\end{lemma}
\pf
The first assertion follows from Proposition \ref{pro:tensor-group}, by taking
\begin{align*}
x&= b v_1\ot w_2 + w_1\ot v_2, & y&=- b v_1\ot w_2 + w_1\ot v_2,\\
v&=w_1\ot w_2, &  w&=b^2 v_1\ot v_2.
\end{align*}
We compute the action of $a_0$ on each one of these generators:
\begin{align*}
	a_0\cdot x
	&= (a-c) \underbrace{b^2 v_1\ot v_2}_{w} - b(a-c) \underbrace{w_1\ot w_2}_{v} 
	\\
	&\hspace*{.8cm}
	-(a + bc) \underbrace{(-bv_1\ot w_2 + w_1\ot v_2)}_{y}\\
	a_0\cdot y
	&=(a-c)\underbrace{b^2 v_1\ot v_2}_{w} + b(a-c) \underbrace{w_1\ot w_2}_{v}+ (bc-a) \underbrace{(bv_1\ot w_2 + w_1\ot v_2)}_{x}\\
	a_0\cdot v
	&=-c \underbrace{b^2v_1\ot v_2}_{w} -b(a+c) \underbrace{bv_1\ot w_2}_{\frac{x-y}{2}} -a \underbrace{w_1\ot w_2}_{v}\\
	a_0\cdot w
	&=a \underbrace{b^2v_1\ot v_2}_{w} + b^2(a+c) \underbrace{w_1\ot v_2}_{\frac{x+y}{2}} + b^2c \underbrace{w_1\ot w_2}_{v}.
\end{align*}
The lemma follows.
\epf

Next, we look for extensions in $\Ext^1(L_j(a,c),L_j(a,c)^{\ot 2})$.

\begin{proposition}\label{pro:Pj}
There is a unique (up to isomorphism) indecomposable module $P_j(a,c)$ so that it fits in an short exact sequence:
\[
0\longrightarrow L_j(a,c)\ot L_j(a,c) \longhookrightarrow P_j(a,c) \longtwoheadrightarrow L_j(a,c) \longrightarrow 0.
\]
It follows that $P_j(a,c)_{|}\simeq S_j^+\oplus S_j^-\oplus M_j^2$ and the action is determined by 
\begin{align}\label{eqn:action-ac}
	[a_0]=\begin{cases}
\begin{psmallmatrix}
	0&-a(b-1)&-ba&b^2a&b& -b^2\\
	-a(1+b)&0&ba&b^2a&0& 0\\
	0&0&-a&b^2a&1& -b^2\\
	0&0&-a&a&1&-1\\
	0&0&0&0&a&-b^2a\\
	0&0&0&0&a&-a
\end{psmallmatrix}, & a=c,\\ 
\\
\begin{psmallmatrix}
	0&-a(1+b)&0&0&1& -b\\
	-a(1-b)&0&0&0&1& b\\
	-2ba&2ba&-a&-b^2a&0& 2b^2\\
	2a&2a&a&a&-2&0\\
	0&0&0&0&a&b^2a\\
	0&0&0&0&-a&-a
\end{psmallmatrix}, & a=-c.
	\end{cases}
\end{align}
In particular, $\dim \Ext^1(L_j(a,c),L_j(a,c)\ot L_0(a,c))=1$.
\end{proposition}
\begin{remark}
We point out that the matrices in \eqref{eqn:action-ac} are written in a ordered basis $\{x,y,v,w,v',w'\}$ to highlight the existence of the submodule $\lg x,y,v,w\rg\simeq L_j(a,c)\ot L_0(a,c)$ and the corresponding quotient $\lg \overline{v'},\overline{w'}\rg\simeq L_j(a,c)$. In the setting of \S\ref{subsec:verymixed}  we  order such a basis as $\{x,y,v,v',w,w'\}$.
\end{remark}
\pf
Let $E=E_j(a,c)\in \Ext^1(L_j(a,c),L)$, with basis $\{x,y,v,w,v',w'\}$ (so $\lg \overline{v'},\overline{w'}\rg=L_j(a,c)$).
There are  $q_1, q_2, q_3, q_4\in\k$ so that
\begin{align*}
	a_0&\cdot v'=q_1 x+ q_2 y +q_3 v +q_4 w + a v' + c w'\\
	&(\Rightarrow a_0\cdot w'=-bq_1 x+ bq_2 y - b^2q_4 v -q_3 w + -b^2c v' -a w').
\end{align*}
That is
\begin{align*}
	[a_0]=\begin{psmallmatrix}
		0&(bc-a)&-\frac{b(a+c)}{2}&\frac{b^2(a+c)}{2}&q_1& -bq_1\\
		-(a+bc)&0&\frac{b(a+c)}{2}&\frac{b^2(a+c)}{2}&q_2& bq_2\\
		-b(a-c)&b(a-c)&-a&b^2c&q_3& -b^2q_4\\
		(a-c)&(a-c)&-c&a&q_4&-q_3\\
		0&0&0&0&a&-b^2c\\
		0&0&0&0&c&-a\\
	\end{psmallmatrix}.
\end{align*}
In terms of Proposition \ref{pro:dimension6}, $\alpha=bc-a$, $\beta=-(a+bc)$, $\delta^+=\delta_-=(a-c \, 0)$,
\begin{align*}
	\eta&=\begin{pmatrix}
		-\frac{b(a+c)}{2}\\q_1
	\end{pmatrix}, & \sigma&=\begin{pmatrix}
		\frac{b(a+c)}{2}\\q_2
	\end{pmatrix} &
	\theta&=\begin{pmatrix}
		-a&0\\q_3 &a
	\end{pmatrix}, & \tau&=\begin{pmatrix}
		-c&0\\q_4&c
	\end{pmatrix}.
\end{align*}

Now, equations \eqref{eqn:general1} hold (recall $f_\pm=\frac{\lambda}{3}$ and $g=\frac{\lambda}{3}(1-b^2)$) and by \eqref{eqn:general2}:
\begin{align*}
	(a-c)(q_1-q_2)&=0, & (a-c)(q_1+q_2)&=2(cq_3-aq_4).
\end{align*}
On the other hand, \eqref{eqn:general3} give, for $\chi=\frac{(a+c)}{2}$:
\begin{align*}
	(bc-a)(q_1-q_2)&=	b(bq_4-q_3)\chi, & (a+bc)(q_1-q_2)& = b(bq_4+q_3)\chi.
\end{align*}
We distinguish two sets of solutions, according to $c=\pm a$:
\begin{itemize}[leftmargin=*]
	\item If $a=c$ (so $\chi=a$) we have a solution for each $q,r\in\k$ as:
	\[
	q_3=q_4=q, \qquad q_2=r, \qquad q_1=bq+r
	\]
	and thus
	\begin{align*}
		[a_0]=\begin{psmallmatrix}
			0&a(b-1)&-ba&b^2a&bq+r& -b(bq+r)\\
			-a(1+b)&0&ba&b^2a&r& br\\
			0&0&-a&b^2a&q& -b^2q\\
			0&0&-a&a&q&-q\\
			0&0&0&0&a&-b^2a\\
			0&0&0&0&a&-a\\
		\end{psmallmatrix}.
	\end{align*}
	\item If $a=-c$ (so $\chi=0$) we have a solution for each $q,r\in\k$ as:
	\[
	q_1=q_2=q, \qquad q_3=r, \qquad q_4=-2q-r.
	\]
	and thus
	\begin{align*}
		[a_0]=\begin{psmallmatrix}
			0&-a(1+b)&0&0&q& -bq\\
			-a(1-b)&0&0&0&q& bq\\
			-2ba&2ba&-a&-b^2a&r& b^2(2q+r)\\
			2a&2a&a&a&-(2q+r)&-r\\
			0&0&0&0&a&b^2a\\
			0&0&0&0&-a&-a\\
		\end{psmallmatrix}.
	\end{align*}
\end{itemize}

It follows that for each $q,r\in\k$ these matrices do define a structure of an $\mA$-module $Q_j(a,c,q,r)$, namely that \eqref{eqn:rels_def}. This can be easily checked by hand, or by using \cite{GAP}.
Next, we observe the following:
\begin{enumerate}
	\item[(a)] If $a=c$, then $Q_j(a,c,q,r)\simeq Q_j(a,c,q+\frac{2r}{b},0)$. 
	\item[(b)] If $a=-c$, then $Q_j(a,c,q,r)\simeq Q_j(a,c,q,0)$. 
\end{enumerate}
Indeed, for case (a), if we set $u=v'-\frac{r}{ba} v$, $z=w'-\frac{r}{ba} w$, we have
\begin{align*}
	a_0\cdot u
	&= (bq+2r)x+(q+\frac{2r}{b})v+(q+\frac{2r}{b})w+au+az.
\end{align*}
Therefore, $Q_j(a,c,q,r)\simeq Q_j(a,c,q+\frac{2r}{b},0)$.
Case (b) follows by setting $u=v'+\frac{r}{a} v$, $z=w'+\frac{r}{a} w$.
Therefore, $Q_j(a,c,q,r)\simeq Q_j(a,c,q,0)$ and  $\dim \Ext^1=1$ (generated by $P_j(a,c)\coloneqq Q_j(a,c,1,0)$).
\epf

In the next section we complete the proof of the claim stating that this module $P_j(a,c)$ is the projective cover of $L_j(a,c)$.

\subsection{The projective covers}

Recall that when $\mu=\c_{\pm j}\frac{\lambda}{3}$ for some (necessarily unique) $j\in\J$, then there is a single six-dimensional simple $\mA$-module $N_j=N_j^{(1)}$ with $N_j=N_j[j]$. 
Otherwise, there are, up to isomorphism two six-dimensional simple $\mA$-modules $N_j^{(1)}$ and $N_j^{(2)}$ for each $j\in\J$.

Recall the definition of the modules $P_j(a,c)$ from Proposition \ref{pro:Pj}.
\begin{theorem}\label{thm:projective}
Let $L$ be an irreducible $\mA_{\lambda,\mu}$-module. Then the projective cover $P(L)$ of $L$ is as follows.
\begin{enumerate}[leftmargin=*]
		\item $\mA(S_0^\pm)$ is the projective cover of the $\A_{\lambda,\mu}$-module $L_0^\pm$.
		\item If $\lambda\neq 0$, then $P_0(a,c)$ is the projective cover of the $\mA$-module $L_0(a,c)$. 
		If $\lambda=0$, then the projective cover of $L_0$ is $\mA(M_0)$.
		\item Assume $\mu=\frac{\lambda}{3}$, then $\mA(S_j^+)$ is the projective cover of $T_j$, $j\geq 1$.		
		\item Assume $\mu=\frac{\lambda}{3}$, then $P_j(a,c)$ is the projective cover of  $L_j(a,c)$. 
		\item Assume $\mu\neq\frac{\lambda}{3}$. If $\mu\neq\c_{\pm j}\frac{\lambda}{3}$, then the six-dimensional modules $N_j^{(1)}$ and $N_j^{(2)}$ are projective. 
		If $\mu=\c_{\pm j}\frac{\lambda}{3}$, then the projective cover of the single six-dimensional module $N_j^{(1)}$ is the induced module $\mA(S_j^+)$.
\end{enumerate}
\end{theorem}
In particular, Items {\it (2)} and {\it (4)} conclude the proof of the claim in \S\ref{sec:claim}.
\pf
{\it (1)} is Proposition \ref{pro:projective0}. For 
{\it (2)}, if $\lambda\neq 0$, then formula
\begin{align*}
72=\dim \mA[0]=\dim L_0^+\dim P(L_0^+)&+\dim L_0^-\dim P(L_0^-)\\
&+\sum_{a,c}\dim L_0(a,c)\dim P(L_0(a,c))
\end{align*}
implies $\dim P(L_0(a,c))=6$ for every $a,c$. As $P_0(a,c)$ is a 6-dimensional projective module, which is, as well, indecomposable and projects onto $L_0(a,c)$, it follows that it is indeed the projective cover. 
On the other hand, if $\lambda=0$, the formula above now gives $\dim P(M_0(0,0))=24$. As this module is necessarily a direct summand of the projective module $\mA_{0,\mu}(M_0)$, the second part of {\it (2)} follows.


{\it (3)} First, we claim that $\mA(S_j^+)\twoheadrightarrow T_j$. Notice that we cannot have $\mA(S_j^+)\twoheadrightarrow L_j(a,c)$ as this would determine four indecomposable summands $P(L_j(a,c))$, each one projecting on each $L_j(a,c)$. As $\dim \mA(S_j^+)=12$, then $\dim P(L_j(a,c))\leq 3$. Now, $L_j(a,c)$ is not projective, which implies $\dim P(L_j(a,c))=3$, but there is not such module. By
\[
72=\dim T_j\dim P(T_j)+\sum_{a,c}\dim L_j(a,c)\dim P(L_j(a,c))
\]
we have that, for each $a,c$ we get that
\[
(\dim P(T_j),\dim P(L_j(a,c)))\in\{(12,6),(8,7),(4,8)\}.
\] 
Now, we also have $\mA(M_j)=\mA_{\lambda,\mu}\ot M_j\twoheadrightarrow L_j(a,c)$, for every $a,c$, and thus $\dim P(L_j(a,c))\geq 6$. Thus $\dim P(T_j)=12$ and therefore $P(T_j)=\mA(S_j^+)$.

Finally, as $\dim P(L_j(a,c))=6$ this yields $P(L_j(a,c))=P_j(a,c)$ and {\it (4)} follows, while 
{\it (5)} is  Corollary \ref{cor:projective-Nj}.
\epf

\section{Extensions and Representation type}\label{sec:extensionsdiagrams}

In this section we compute the spaces $\Ext^1(L,L')$ between two simple $\mA_{\lambda,\mu}$-modules $L$ and $L'$. That is, we compute all $\mA_{\lambda,\mu}$-modules $E$ fitting in:
\[
0\to L'\hookrightarrow E\twoheadrightarrow L\to 0.
\]
With this information, we describe the Gabriel quiver of $\mA_{\lambda,\mu}$ and its separated diagram in \S\ref{sec:diagrams}; to determine its representation type. 
We recall once again the notation from  \S\ref{sec:notation}.

\subsection{$\Ext$ computations}\label{sec:extensions}

We begin with a simple observation. To do this, we recall from Proposition \ref{pro:j-unico} the definition of the $j$th component $M[j]=M_{|}[j]$ cf.~\eqref{eqn:j-part} of a given $\A$-module $M$. The following is straightforward.
\begin{lemma}
	Let $S=S[j]$ and $T=T[k]$ two simple $\A$-modules. 
	If $j\neq k$, then $\dim\Ext^1(S,T)=0$.
\end{lemma}
\pf
If $M\in \Ext^1(S,T)$, then $M[j]=S$ and $M[k]=T$, so $M\simeq S\oplus T$.
\epf

In this part we shall look at the general, non-graded, case, where $\lambda$ and $\mu$ are not simultaneously zero. We include in \S\ref{sec:ext-graded} the graded case $\lambda=\mu=0$.

We invoke Lemma \ref{lem:t0} to deal with extensions between one-dimensional modules all at once. We recall the modules $T_0^{\eps\eps'}$ in \eqref{eqn:t0}, for $\eps\neq\eps'\in\{\pm\}$.

\begin{corollary} 
		If $M\in\Ext^1 (L_0^{\eps},L_0^{-\eps})$, then $M\simeq T_0^{\eps\eps'}$. In particular, we have that $\dim\Ext^1(L_0^{\eps},L_0^{\eps'})=\delta_{\eps,-\eps'}$.
	\qed
\end{corollary}

Next, we deal with extensions between one-dimensional modules $L_0^\pm$ and the two-dimensional modules supported on $M_0$. We use Proposition \ref{pro:sumadedim2} to deal with extensions between modules of type $L_j(a,c)$. Recall that in this last case we allow $j\in\J$ when $\mu=\frac{\lambda}{3}$.

\begin{lemma}\label{lem:extSM} The following holds.
	\begin{enumerate}[leftmargin=*]
		\item $\dim \Ext^1(L_0^\eps,L_0(a,c))=\delta_{c,-\eps a}$.
		If $M\in \Ext^1(L_0^\eps,L_0(a,\eps a))$, then there is $\alpha\in\k$ such that $[a_0]=\begin{psmallmatrix}
			a&\eps a&\alpha\\
			-\eps a&-a&-\eps\alpha\\
			0&0&0.
		\end{psmallmatrix}$, $\alpha\in\k$.
		\item $\dim \Ext^1(L_0(a,c), L_0^\eps)=\delta_{c,\eps a}$.
		If $M\in \Ext^1(L_0(a,\eps a), L_0^\eps)$, then there is $\alpha\in\k$ such that   $[a_0]=\begin{psmallmatrix}
			0&\alpha&-\eps\alpha\\
			0&a&-\eps a\\
			0&\eps a&-a
		\end{psmallmatrix}$, $\alpha\in\k$.
		\item[(3i)]If $\lambda\neq 0$, then $\dim\Ext^1(L_j(a,c), L_j(a',c'))=0$, for every selection $a,c,a',c'\in\{\pm\sqrt{\lambda/3}\}$.
\item[(3ii)] If $\lambda=0$, then 
$\dim\Ext^1(L_0, L_0)=2$. If $M\in\Ext^1(L_0, L_0)$ then $[a_0]=\begin{psmallmatrix}
	0&0&\alpha&-\beta\\
	0&0&\beta&-\alpha\\
	0&0&0&0\\
	0&0&0&0\\
\end{psmallmatrix}$, for some $\alpha, \beta\in\k$.		
	\end{enumerate}
\end{lemma}
\pf
{\it (1)} and {\it (2)} are easy to check and {\it (3i)} follows from Proposition \ref{pro:sumadedim2}, as there are no indecomposable modules arising from the sum of two $\G_{3,\ell}$-modules of dimension 2. Finally, {\it (3ii)} is clear.
\epf

\begin{remark}
Assume $\lambda=0$. Then Items {\it (1)} and {\it (2)}  from Lemma \ref{lem:extSM} read $\dim \Ext^1(L_0^\eps,L_0)=\dim \Ext^1(L_0, L_0^\eps)=1$. In turn, Item {\it (3ii)} contains Example \ref{exa:extensionM0}, for $\alpha=1, \beta=0$.
\end{remark}

Finally, we turn to the modules $T_j$ and $L_j(a,c)$, $j\in\J$, which correspond exclusively to the case $\mu=\frac\lambda3$.

\begin{lemma}\label{lem:extTM}
Assume $\mu=\frac\lambda3$.
	\begin{enumerate}[leftmargin=*]
			\item $\dim \Ext^1(L_j(a,c), T_j)=1$.
	If $M\in \Ext^1(L_j(a,c), T_j)$, then 
	\begin{align}\label{eqn:a0-extMjTj}
		[a_0]=\left(\begin{smallmatrix*}
			0&\mu(1-\zeta^j)&-\beta(a+cb)&b\beta(a+cb)\\
			1&0&\beta&b\beta\\
			0&0&a&-\zeta^j c\\
			0&0&c&-a
		\end{smallmatrix*}\right), \qquad \beta\in\k.
	\end{align}
	\item $\dim \Ext^1(T_j, L_j(a,c))=
	2 $. If $M\in \Ext^1(T_j, L_j(a,c))$, then 
	\[
	[a_0]=\left(\begin{smallmatrix*}
		a&-\zeta^j c&-b\beta&b\gamma\\
		c& -a&\beta&\gamma\\
		0&0&0&\mu(1-\zeta^j)\\
		0&0&1&0
	\end{smallmatrix*}\right), \qquad \beta, \gamma\in \k.
	\]
	\item $\dim\Ext^1(T_j,T_j)=0$.
	\end{enumerate}
\end{lemma}

\pf
We deal with \textit{(1)}, as \textit{(2)} follows similarly. Set $M=\lg x, y, v, w\rg$ so that
$$ 0\longrightarrow {} T_k \longhookrightarrow M \longtwoheadrightarrow L_j(a,c) \longrightarrow 0,$$
that is $\lg x,y\rg\simeq T_j$ and $\lg v,w\rg/T_j\simeq L_j(a, c)$.
That is, we have that $a_0\cdot x=y$, $a_0\cdot y=\mu(1-\zeta^j)x$ and there are $\alpha,\beta\in\k$ so that
\begin{align*}
	a_0\cdot v&=\alpha x+\beta y+av+cw.
\end{align*}
In particular, 	$a_0\cdot w=-\alpha b x + \beta b y -\zeta^j v- aw$.

Applying the relation $a_0^2=\mu(1-t^2)$ on $v,w$ we get that:
\begin{align*}
		(a-cb)\alpha+\mu(1-\zeta^j)\beta&=0, & 
		\alpha + (a+cb)\beta&=0. 
\end{align*}
Thus we get that $\alpha=-\beta(a+cb)$ and thus $[a_0]$ is as in \eqref{eqn:a0-extMjTj}.
It follows that for this matrix the relation $a_0a_1+a_1a_2+a_2a_0=\lambda(1-st^2)$ is satisfied. Hence this defines an $\mA_{\lambda,\mu}$-module for each $\beta$. It is clear that isomorphism classes are classified by $\beta=0$ (direct sum) and $\beta=1$.

As for \textit{(3)}, let  $M=\lg x, y, x', y'\rg\in \Ext^1(T_j, T_j)$ such that for some $\alpha,\beta\in\k$:
\begin{align*}
	a_0\cdot x &= y, & a_0\cdot y&=\mu(1-\zeta^j)x, &	a_0\cdot x' &=\alpha x + \beta y+ y'.
\end{align*}
Hence $a_0\cdot y'=a_0\cdot (a_0\cdot x'-\alpha x - \beta y) =-\beta\mu(1-\zeta^j) x - \alpha y + \mu(1-\zeta^j)x'$.

Setting $x''=x'-\beta x$ and  $y''=a_0\cdot x''=\alpha x + y'$ we have
\begin{align*}
	a_0\cdot x''&=y'', && a_0\cdot y''=\mu(1-\zeta^j)x''.
\end{align*}
and thus $M\simeq T_j\oplus T_j$.
\epf

We are left with the case of the 6-dimensional simple modules $N_j$, $j\in\J$, for $\mu\neq\frac{\lambda}{3}$. Notice that we are only interested in the case in which there is a unique such module for each $j$, as otherwise these modules are projective (and hence have no non-trivial extensions).

Hence we assume next that $\mu=\c_{\pm j}\frac{\lambda}{3}$ holds, that is there a single 6-dimensional module $N_j^{(1)}$. 
In particular, notice that in this setting:
	\begin{align}\label{eqn:c1c2}
		c_1t_1-c_2t_2=0.
	\end{align}
	Indeed, this is clear if $c_1=c_2=0$. If, otherwise, $c_1=1, c_2=g$ and $t_1=h/2$, $t_2=2f_+/h$ (with $h^2=4f_+g$), this identity is satisfied as well.

Now we can compute the extensions.
\begin{lemma}
Assume $\mu=\c_{\pm j}\frac{\lambda}{3}$. Then $\dim\Ext^1(N_j^{(1)}, N_j^{(1)})=2$.
	%
\end{lemma}
\pf
Consider two copies of $N_j$, with bases 
\[
\{x,y,v_1,v_2,w_1,w_2\}, \quad \{\tilde{x},\tilde{y},\tilde{v}_1,\tilde{v}_2,\tilde{w}_1,\tilde{w}_2\}
\]
respectively. Let $M\in\dim\Ext^1(N_j, N_j)$ with  $[a_0]$ as in  \eqref{eqn:a0matrix}
for the basis
\[
B_N=\{x,\tilde{x},y,\tilde{y},v_1,v_2,\tilde{v}_1,\tilde{v}_2,w_1,w_2,\tilde{w}_1,\tilde{w}_2\}.
\] 
Hence there are scalars $\alpha_1$, $\beta_1$, $\delta_1$, $\delta_2$, $\delta_3$, $\delta_4$, $\eta_1$, $\eta_2$, $\sigma_1$, $\sigma_2$, $\theta_1$, $\theta_2$, $\theta_3$, $\theta_4$, $\tau_1$, $\tau_2$, $\tau_3$, $\tau_4$ so that 
\begin{align*}
	\alpha=&\left[\begin{smallmatrix*} c_2 & 0\\ \alpha_1 & c_2 \end{smallmatrix*}\right], & \beta&=\left[\begin{smallmatrix*} c_1 & 0\\ \beta_1 & c_1 \end{smallmatrix*}\right], &
	\delta_+=&\left[\begin{smallmatrix*} 1&0&0&0\\ \delta_1 & \delta_2&1&0 \end{smallmatrix*}\right], &
	\delta_-=&\left[\begin{smallmatrix*} 0&1&0&0\\ \delta_3 & \delta_4&0&1 \end{smallmatrix*}\right]\\
	\eta=& \left[\begin{smallmatrix*} f&0\\0&0\\\eta_1&f\\\eta_2&0 \end{smallmatrix*}\right], & 
	\sigma=& \left[\begin{smallmatrix*} 0&0\\-f&0\\\sigma_1&0\\\sigma_2&-f\\\end{smallmatrix*}\right], & 
	\theta=& \left[\begin{smallmatrix*} 0&t_2&0&0\\ t_1&0&0&0\\ \theta_1&\theta_2&0&t_2\\\theta_3&\theta_4&t_1&0  \end{smallmatrix*}\right], & 
	\tau=& \left[\begin{smallmatrix*} 0&t_2'&0&0\\ t_1'&0&0&0\\ \tau_1&\tau_2&0&t_2'\\ \tau_3&\tau_4&t_1'&0\end{smallmatrix*}\right]
\end{align*}
and the scalars $c_1, c_2, t_1, t_2$ -and $t_1'=\frac{t_1-c_2}{b}$, $t_2'=\frac{c_1-t_2}{b}$- are determined by each $N_j$ structure.
From equations \eqref{eqn:general1}--\eqref{eqn:general3} we obtain the following restrictions.
By  \eqref{eqn:general1}, we get
\begin{align*}
\alpha\beta&=\beta\alpha=g\id & \implies	&\begin{cases*}
		\alpha_1c_1=-\beta_1c_2,	
	\end{cases*}
	\\
\theta^2&=f_+\id &\implies
	&\begin{cases*}
		t_1\theta_2=-t_2\theta_3, \\ t_2(\theta_1+\theta_4)=0, \\ t_1(\theta_1+\theta_4)=0,
	\end{cases*}
	\stackrel{t_1t_2\neq0}{\Longleftrightarrow}
	\begin{cases*}
		\theta_4=-\theta_1, \\ \theta_3=-\frac{t_1}{t_2}\theta_2,
	\end{cases*}
	\\
\tau^2&=f_-\id &\implies
	&\begin{cases*}
		\tau_2(t_1-c_2)=-(c_1-t_2)\tau_3, \\ \tau_1(c_1-t_2)=-(c_1-t_2)\tau_4, \\ \tau_4(t_1-c_2)=-(t_1-c_2)\tau_1,
	\end{cases*}
	\stackrel{c_1\neq t_2}{\Longleftrightarrow}
	\begin{cases*}
		\tau_4=-\tau_1, \\ \tau_3=-\frac{t_1-c_2}{c_1-t_2}\tau_2,
	\end{cases*}
\end{align*}
We look now at \eqref{eqn:general2}, which yields:
\begin{align*}
	\sigma_1&=\delta_2f, & \sigma_2&=\delta_4f , & \eta_1&=-\delta_1f, & \eta_2&=-\delta_3f.
\end{align*}
As well, $\eta\delta_++\sigma\delta_-=\tau\theta-\theta\tau$ gives
\begin{align*}	
	\begin{cases*}
		\frac{c_1}{t_2(c_1-t_2)}\underbrace{(c_1t_1-c_2t_2)}_{=0}(\tau_2-\frac{\theta_2}{b})=0,\\
		\delta_2=\frac{1}{f}\left(t_2\tau_1-\frac{c_1-t_2}{b}\theta_1\right), \qquad 
		\delta_3=\frac{1}{f}\left(t_1\tau_1+\frac{c_2-t_1}{b}\theta_1\right).
	\end{cases*}
\end{align*}
Finally, $\beta\delta_--\delta_+(b\tau+\theta)=0=\alpha\delta_++\delta_-(b\tau-\theta)$ from \eqref{eqn:general3}, implies
\begin{align*}	
	\begin{cases*}
		\underbrace{(bf-h+2f_+)}_{=0}(\theta_1-b\tau_1)=0\\
		\beta_1=\theta_2+b\tau_2+c_1(\delta_1-\delta_4)\\
		\alpha_1=-\frac{t_1}{t_2}\theta_2+\frac{t_1-c_2}{c_1-t_2}b\tau_2-c_2(\delta_1-\delta_4)
	\end{cases*}
\end{align*}
where we check that these expressions for $\alpha_1$ and $\beta_1$ are consistent with the previous ones. Similarly, the remaining equations from \eqref{eqn:general3} are also satisfied.
Hence we get to
\begin{align*}
	\alpha&=\left[\begin{smallmatrix*} c_2 & 0\\ -\frac{t_1}{t_2}\theta_2+\frac{t_1'}{t_2'}b\tau_2-c_2(\delta_1-\delta_4) & c_2 \end{smallmatrix*}\right], &
	\beta&=\left[\begin{smallmatrix*} c_1 & 0\\ \theta_2+b\tau_2+c_1(\delta_1-\delta_4) & c_1 \end{smallmatrix*}\right],
\end{align*}
\begin{align*}
	\delta_+&=\left[\begin{smallmatrix*} 1&0&0&0\\ \delta_1 & \frac{t_2}{f}\tau_1-\frac{t_2'}{f}\theta_1&1&0 \end{smallmatrix*}\right],& 
	\eta&= \left[\begin{smallmatrix*} f&0\\0&0\\-\delta_1f&f\\-t_1\tau_1+t_1'\theta_1&0 \end{smallmatrix*}\right],  &
	\theta=& \left[\begin{smallmatrix*} 0&t_2&0&0\\ t_1&0&0&0\\ \theta_1&\theta_2&0&t_2\\-\frac{t_1}{t_2}\theta_2&-\theta_1&t_1&0  \end{smallmatrix*}\right], \\
		\delta_-&=\left[\begin{smallmatrix*} 0&1&0&0\\ \frac{t_1}{f}\tau_1+\frac{t_1'}{f}\theta_1 & \delta_4&0&1 \end{smallmatrix*}\right], &
	\sigma&= \left[\begin{smallmatrix*} 0&0\\-f&0\\ t_2\tau_1-t_2'\theta_1&0\\\delta_4f&-f\\\end{smallmatrix*}\right], &
	\tau=& \left[\begin{smallmatrix*} 0&t_2'&0&0\\ t_1'&0&0&0\\ \tau_1&\tau_2&0&t_2'\\ -\frac{t_1'}{t_2'}\tau_2&-\tau_1&t_1'&0\end{smallmatrix*}\right].
\end{align*}

Now we perform the following base change, from $B_N$ to 
\[
\overline{B}_N=\{x,\bar{x},y,\bar{y},v_1,v_2,\bar{v}_1,\bar{v}_2,w_1,w_2,\bar{w}_1,\bar{w}_2\}.
\]
where
\begin{align*}
	\bar{x}&=\tilde{x}-\delta_1 x, & \bar{y}&=\tilde{y}-\delta_4 y, & \bar{v}_1&=\tilde{v}_1-\frac{\theta_1}{t_1}v_2, \\ 
	\bar{v}_2&=\tilde{v}_2+\frac{\theta_2}{t_2}v_2, & \bar{w}_1&=\tilde{w}_1-\frac{\theta_1}{t_1}w_2, & \bar{w}_2&=\tilde{w}_2+\frac{\theta_2}{t_2}w_2.
\end{align*}
It is easy to check that, in this new basis, the matrices only depend on two parameters $\tau_1$ and $\tau_2$ and thus the matrix of $a_0$ becomes
\begin{align*}
	[a_0]=\begin{psmallmatrix}
		0&0&c_2&\frac{t_1'}{t_2'}b\tau_2&f&0&0&-t_1\tau_1&-bf&0&0&bt_1\tau_1\\ 
		0&0&0&c_2&0&0&f&0&0&0&-bf&0\\
		c_1&b\tau_2&0&0&0&-f&t_2\tau_1&0&0&-bf&bt_2\tau_1&0\\
		0&c_1&0&0&0&0&0&-f&0&0&0&-bf\\
		-b&0&0&\frac{t_1}{f}b\tau_1&0&t_1&0&0&0&-b^2 t_1'&-b^2\tau_1&\frac{t_1'}{t_2'}b^2\tau_2\\
		0&-\frac{t_2}{f}b\tau_1&b&0&t_2&0&0&0&-b^2 t_2'&0&-b^2\tau_2&b^2\tau_1\\
		0&-b&0&0&0&0&0&t_1&0&0&0&-b^2t_1'\\
		0&0&0&b&0&0&t_2&0&0&0&-b^2 t_2'&0\\
		1&0&0&\frac{t_1}{f}\tau_1&0&t_1'&\tau_1&-\frac{t_1'}{t_2'}\tau_2&0&-t_1&0&0\\
		0&\frac{t_2}{f}\tau_1&1&0&t_2'&0&\tau_2&-\tau_1&-t_2&0&0&0\\
		0&1&0&0&0&0&0&t_1'&0&0&0&-t_1\\
		0&0&0&1&0&0&t_2'&0&0&0&-t_2&0
	\end{psmallmatrix}.
\end{align*}
Thus, we obtain an extension $N_j(\tau_1, \tau_2)\in\ext^1(N_j,N_j)$, for each $\tau_1,\tau_2\in \k$. A routine calculation shows $N_j(1,0)\not\simeq N_j(0,1)$.
\epf

\subsubsection{The graded case}\label{sec:ext-graded}
We look into the graded case $\lambda=\mu=0$. 
Recall that the simple $\mA_{0,0}$-modules are in one-to-one correspondence with $\widehat{\G_{3,\ell}}$.
The proof is straightforward.

\begin{lemma}
	If $\lambda=\mu=0$ then for each $j=0, \dots, \ell-1$:
	\begin{enumerate}[leftmargin=*]
		\item $\dim\Ext^1(L_j^\eps, L_j^{\eps'})=\delta_{\eps,-\eps'}$.
		\item $\dim\Ext^1(L_j^\eps,L_j)=\dim\Ext^1(L_j, L_j^\eps)=1$
		\item $\dim\Ext^1(L_j, L_j)=2$. \qed
	\end{enumerate}
\end{lemma}

\subsection{Diagrams}\label{sec:diagrams}

In this part we use the data obtained above to compute the Gabriel quiver and the corresponding Separated Ext-Quiver (SEQ) for $\mA_{\lambda,\mu}$. This allows us to establish the following, see \S\ref{sec:proof-diagrams} for a proof.

\begin{proposition}\label{pro:repr-type}
	The algebra $\mA_{\lambda,\mu}$ is not of finite representation type. If $\lambda=0$ or $\mu=\frac\lambda3$, then it is of wild representation type.\qed
\end{proposition}

To describe the diagrams, we introduce some notation. We set $\hslash=\sqrt\frac\lambda3$. As well, we shall use the following drawings.
\begin{itemize}[leftmargin=*]
	\item[$\oplus$] a collection of four vertices, each one corresponding to a module
	$L_0(\hslash,\hslash)$, $L_0(-\hslash,-\hslash)$, $L_0(-\hslash,\hslash)$ or $L_0(\hslash,-\hslash)$,
	\item[$\overt$] to refer to two vertices together, corresponding to 
	$L_0(\hslash,\hslash)$, $L_0(-\hslash,-\hslash)$, 
	\item[$\ominus$] for two vertices associated to $L_0(\hslash,-\hslash)$, $L_0(-\hslash,\hslash)$.
\end{itemize}

We distinguish three cases (graded, Hochschild, pure), namely:
\begin{enumerate}[leftmargin=*]
	\item $\lambda=\mu=0$, see \S\ref{ext-0}.
	\item $\mu=\frac\lambda3\neq0$, see \S\ref{ext-lambda/3}.
	\item $\mu\neq\frac\lambda3\neq0$. In turn, this splits into two cases, according to wether $\mu=\c_{\pm j}\frac{\lambda}{3}$ for some $j\in\J$. Namely, if there are two non-isomorphic 6-dimensional irreducible modules, see \S\ref{ext-twoNj}, or there is $j\in\J$ for which there is a single one, see \S\ref{ext-oneNj}.
\end{enumerate}

\subsection{$\lambda=\mu=0$.}\label{ext-0}

The Ext-Quiver is given by $\ell$ copies $\bigsqcup_{j=0}^{\ell-1}G_j$, where

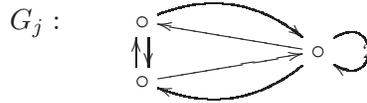
\begin{equation*}
	G_j :\qquad \xymatrix@R=-.1pt{ \circ\ar@<0.5ex>[dd]\ar@/^1pc/[rrd]&& \\
		&& \circ\ar[llu]\ar@/^1pc/[lld]\ar@(ur,dr)[]|{2}\\
		\circ\ar[rru]\ar@<0.5ex>[uu]&&}
\end{equation*}

So, the Separated Ext-Quiver is the disjunct union $\bigsqcup_{j=0}^{\ell-1}F_j$, where
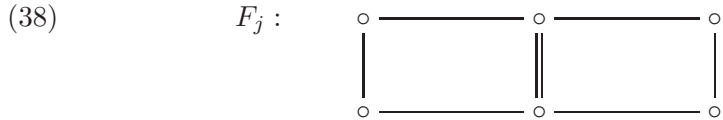
\begin{equation}\label{eqn:diagramFj}
	F_j :\qquad 	\xymatrix{ \circ\ar@{-}[rr]&&\circ\ar@{-}[rr]&&\circ\ar@{-}[d]\\
		\circ\ar@{-}[u]&&\circ\ar@{-}[ll]\ar@{=}[u]&&\circ\ar@{-}[ll]
	}
\end{equation}

\subsection{$\mu=\frac\lambda3\neq0$.}\label{ext-lambda/3}
The Ext-Quiver has $\ell$ connected components, associated to sets of vertices $V_0=\{L_0^\eps,T_0,L_0(a,c):a,c=\eps\hslash\}$ and $V_j=\{T_j,L_j(a,c):a,c=\eps\hslash\}$, $j=1,\dots\ell-1$.
The component  associated to $V_0$ is
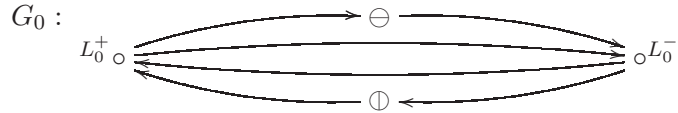
\begin{equation*}
	G_0: \xymatrix@R=-.1pt{&&&\ominus\ar@/^/[drrr]&&&\\
		^{L_0^+}\circ\ar@/^/[rrrrrr]\ar@/^/[urrr]\ar@/^/[urrr]&&&&&&\circ^{L_0^-}\ar@/^/[llllll]\ar@/^/[dlll]\\
		&&&\overt\ar@/^/[ulll]&&&}
\end{equation*}
and thus the corresponding component $E_0$ of the SEQ is $D_5^{(1)}\sqcup D_5^{(1)}$, where 
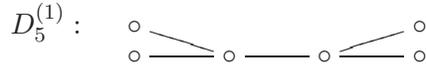
\begin{equation*}
	D_5^{(1)}: \quad\xymatrix@R=-.1pt{ \circ\ar@{-}[rd]	& & &\circ\ar@{-}[ld]\\
		\circ\ar@{-}[r]&\circ\ar@{-}[r]&\circ\ar@{-}[r]&\circ
	}
\end{equation*}
Now, the associated component to each $V_j$, $j\in\J$, is the quiver:
\begin{equation}\label{eq: Gj}
	G_j : \qquad	\xymatrix{ \oplus\ar@/^1pc/[rrr]&&&\circ\ar@/^/[lll]|{2}}
\end{equation}
and thus the corresponding component of SEQ is 
\begin{equation}\label{eqn:diagramDD4}
	E_j : \qquad \xymatrix@R=-.1pt{ \circ\ar@{-}[rd]&&\circ\ar@{-}[ld]\\
		\circ&\circ\ar@{-}[r]\ar@{-}[l]&\circ
	}
	\qquad
	\xymatrix@R=-.1pt{ \circ\ar@{=}[rd]&&\circ\ar@{=}[ld]\\
		\circ&\circ\ar@{=}[r]\ar@{=}[l]&\circ
	}
\end{equation}
Let us call $\mathbb{D}_4^{(1)}$ the diagram in the right of \eqref{eqn:diagramDD4}, so $E_j=D_4^{(1)}\sqcup \mathbb{D}_4^{(1)}$. Then the SEQ is
$(D_5^{(1)})^2\sqcup(D_4^{(1)})^{\ell-1}\sqcup(\mathbb{D}_4^{(1)})^{\ell-1}$.

\subsubsection{}\label{ext-twoNj}
For $j=0$, the component $G_0$ is either the diagram in \S\ref{ext-lambda/3}, if $\lambda\neq 0$ or the diagram in \S\ref{ext-0} if $\lambda=0$.
The components of the Ext-Quiver for each $j\in\J$ are described by the Dynkin diagram $A_1\sqcup A_1$:
\begin{equation*}
	G_j: \qquad \xymatrix@1{
		\circ && \circ}
\end{equation*}
Hence the SEQ is for $\lambda\neq 0$ or $\lambda=0$ are  $(D_5^{(1)})^2\sqcup A_1^{4(\ell-1)}$ or $F_0\sqcup A_1^{4(\ell-1)}$, respectively, for $F_0$ is as in \eqref{eqn:diagramFj}:

\subsubsection{}\label{ext-oneNj}
Here we assume that $\mu=\c_{\pm j}\frac{\lambda}{3}$ for some (unique) $j\in\J$.

Again, for $j=0$, the component $G_0$ is either the diagram in \S\ref{ext-lambda/3}. To describe the components for $j\in\J$, let us set
\begin{align*}
	\mathbb{A}_0^{(1)}&: 
		\xymatrix@1{
			\circ\ar@(ur,dr)[]|{2}}
\end{align*}
In particular, observe that the separated component of this diagram is $B_2$:
\begin{equation*}
	\xymatrix@1{
		\circ\ar@{=}[rr] && \circ}
\end{equation*}
We have that $G_j=\mathbb{A}_0^{(1)}$ and $G_k=	A_1^2$, for $k\neq j$.
Hence the SEQ is, respectively: $
	(D_5^{(1)})^2\sqcup B_2\sqcup A_1^{4(\ell-2)}$ or $(D_5^{(1)})^2\sqcup  B_2^{\ell-1}$.

\subsection{Proof of Proposition \ref{pro:repr-type}}\label{sec:proof-diagrams} 

If $\lambda=0$, then the separated quiver contains a component $F_0$ as in \eqref{eqn:diagramFj}, which is not of finite nor affine type. The same holds when $\mu=\frac{\lambda}{3}$, for a component $\mathbb{D}_4^{(1)}$ as in \eqref{eqn:diagramDD4}. If $\mu\neq\frac\lambda{3}\neq 0$, then the separated quiver contains at least a component of affine type $D_5^{(1)}$ and any other component is of affine (as $B_2$) or finite type (as $A_1$).\qed

\end{document}